\documentstyle[amstex,amssymb,12pt]{amsart}
\numberwithin{equation}{section}
\newtheorem{dfn}{Definition}[section]
\newtheorem{prop}[dfn]{Proposition}
\newtheorem{thm}[dfn]{Theorem}
\newtheorem{lem}[dfn]{Lemma}

\newtheorem{rem}[dfn]{Remark}

\begin{document}
\title{The multiple gamma functions and the multiple $q$-gamma
functions}
\author{Kimio Ueno and Michitomo Nishizawa}
\address{Department of Mathematics,
 School of Science and Engineering,
 Waseda University}
\email{uenoki@@cfi.waseda.ac.jp, 694m5035@@cfi.waseda.ac.jp}
\maketitle

\footnote[0]{1995 {\it Mathematical Subject Classification.} 33B15, 33D05.}

\begin{abstract}
We give an asymptotic expansion ({\it the higher Stirling formula})
and an infinite product representation ({\it the Weierstrass product
representation}) of the Vign\'{e}ras multiple gamma functions by considering
the classical limit of the multiple $q$-gamma functions.
\end{abstract}

\section{Introduction}
The multiple gamma function was introduced by Barnes. It is defined
to be an infinite product regularized by the multiple Hurwitz zeta functions
\cite{bar1}, \cite{bar2}, \cite{bar3}, \cite{bar4}. After his discovery,
many mathematicians have studied this  function:
Hardy \cite{har1}, \cite{har2} studied
this function from his viewpoint of the theory of elliptic functions,
and Shintani \cite{shi1},\cite{shi2} applied it to the study on the 
Kronecker limit formula for zeta functions attached to certain 
algebraic fields.\par
 In the end of 70's, Vign\'{e}ras \cite{vig} redefined the multiple
gamma function to be a function satisfying the generalized Bohr-Morellup
theorem, Furthermore, Vign\'{e}ras \cite{vig}, Voros \cite{vor},
Vardi \cite{var} and  Kurokawa
\cite{kur1}, \cite{kur2}, \cite{kur3}, \cite{kur4} showed that it plays
an essential role to express gamma factors of the  Selberg zeta functions
of compact Riemann surfaces and the determinants of the Laplacians on some 
Riemannian manifolds.\par
As we can see from these studies, the multiple gamma functions are
fundamental for the analytic number theory: See also \cite{kur5},
\cite{man}. However we do not think that the theory of the multiple gamma
functions has been fully explored.\par
 On the other hand, the second author of this paper introduced 
a $q$-analogue of the Vign\'{e}ras multiple gamma functions and showed it to be
characterized by a $q$-analogue of the generalized Bohr-Morellup
theorem \cite{vig}.\par
 In this paper, we will establish an asymptotic expansion formula
({\it the higher Stirling formula}) and an infinite product representation
({\it the Weierstrass product representation}) of the Vign\'{e}ras multiple 
gamma functions by considering the classical limit of the multiple
$q$-gamma functions. In order to get these results, we will use the method
developed in \cite{un}. Namely, by making use of the Euler-MacLaurin summation
formula, we derive the Euler-MacLaurin expansion of the multiple $q$-gamma
functions. Taking the classical limit, we are led to the
Euler-MacLaurin expansion of the Vig\'{n}eras multiple gamma functions.
The higher Stirling formula and the Weierstrass product representation
follow from this expansion formula \par
 This paper is organized as follows. In Section 2, we give a survey of 
the multiple gamma functions and its $q$-analogue. In Section
3, we derive the Euler-MacLaurin expansion of the multiple $q$-gamma functions
by using the Euler-MacLaurin summation formula. In Section 4, we consider
the classical limit of the multiple $q$-gamma functions rigorously and give
an asymptotic expansion formula of the Vign\'{e}ras multiple gamma functions. 
In Section 5, the Weierstrass product representation of this function is 
derived.\par 

\vspace{12pt}

\noindent{\bf Acknowledgement.}
The first author is partially supported by Grant-in-Aid for Scientific
Research on Priority Area 231 ''Infinite Analysis'' and by Waseda
University Grant for Special Research Project 95A-257.



\section{A survey of the multiple gamma function and the multiple
$q$-gamma function}
\subsection{The gamma function}
The following are well-known facts in the classical analysis:
The Bohr-Morellup theorem says that the gamma function $\Gamma(z)$
is characterized by the three conditions,
   \begin{align*}
      & \mbox{(1)} \quad \Gamma(z+1)=z \Gamma(z),\\
      & \mbox{(2)} \quad \Gamma(1)=1,\\
      & \mbox{(3)} \quad \frac{d^2}{dz^2}\log\Gamma(z+1)\geq0
        \quad \mbox{for} \quad z\geq0.
   \end{align*}
The gamma function is meromorphic on {\bf C}, and has an infinite product
representation
   \begin{equation}
     \Gamma(z+1)=e^{-\gamma x}
       \prod_{n=1}^{\infty}
        \left\{
        \left(1+\frac{z}{n}\right)^{-1}
        e^{\frac{z}{n}}
        \right\},
   \end{equation}
where $\gamma$ is the Euler constant. This is usually called the Weierstrass
product formula.\par
Another representation of the gamma function is derived from the Hurwitz
zeta function:
  \begin{equation}
    \zeta(s,z):= \sum_{k=0}^{\infty}
      \frac{1}{(z+k)^{s}} \quad \mbox{for} \quad \Re s > 1.
  \label{eqn:zgam}\end{equation} 
It is well-known that
  \begin{equation*}
    \frac{\Gamma(z)}{\sqrt{2\pi}}=\exp(\zeta'(0,z)),
  \end{equation*}
where $\zeta'(0,z)= \frac{d}{ds}\zeta(s,z)|_{s=0}.$\par
The gamma function has an asymptotic expansion formula, {\it i.e.
the Stirling formula,}
  \begin{eqnarray*}
    & &\log \Gamma(z+1) \sim
      \left(z+\frac{1}{2}\right)\log(z+1) - (z+1) - \zeta'(0)\\
    & &\qquad + \sum_{r=1}^{\infty} \frac{B_{2r}}{[2r]_{2}}
    \frac{1}{(z+1)^{2r-1}},
  \end{eqnarray*}
as $z\to\infty$ in a sector $\Delta_{\delta}:= \{z\in{\bf C}|
|\arg z|<\pi-\delta\}$ $(0<\delta<\pi),$
where
     $$\frac{z e^{tz}}{e^{z}-1}
        = \sum_{n=0}^{\infty} \frac{B_{n}(t)}{n!} z^{n},$$
$B_{k}=B_{k}(0)\quad(\mbox{the Bernoulli number})$,  $\zeta(s)$
is the Riemann zeta function, $\zeta'(s)= \frac{d}{ds}\zeta(s)$ and
$[x]_{r}=x(x-1)\cdots (x-r+1)$. Note that $\zeta'(0)=-\log\sqrt{2\pi}$.
\par


\subsection{The Barnes $G$-function}
Barnes \cite{bar1} introduced  the function $G(z)$ which satisfies
   \begin{align*}
      & \mbox{(1)} \quad G(z+1)= \Gamma(z) G(z),\\
      & \mbox{(2)} \quad G(1)=1,\\
      & \mbox{(3)} \quad \frac{d^3}{dz^3}\log G(z+1)\geq0
        \quad \mbox{for} \quad z\geq0,
   \end{align*}
and he called this  the ``$G$-function''. He proved that the $G$-function 
has an infinite product representation.
  $$ G(z+1)
      = e^{-z\zeta'(0)-\frac{z^{2}}{2}\gamma
      - \frac{z^{2}+z}{2}}
     \prod_{k=1}^{\infty} \left\{\left(
        1+\frac{z}{k}\right)^{k}
        \exp\left(-z+\frac{z^{2}}{2k}\right)\right\}$$
and an asymptotic expansion
  \begin{equation}
    \log G(z+1)\sim\left(\frac{z^{2}}{2}-\frac{1}{12}\right)\log(z+1)
      -\frac{3}{4}z^{2}-\frac{z}{2}+\frac{1}{3}+z\zeta'(0)
      -\log A + O(\frac{1}{z})
  \label{eqn:asG}\end{equation}
as $z\to\infty$ in the sector $\Delta_{\delta}$, where $A$ is called the
Kinkelin constant. Voros showed this constant can be written with the first
derivative of the Riemann zeta function
(cf \cite{vor}, \cite{var})
  $$ \log A = -\zeta'(-1)+\frac{1}{12}.$$


\subsection{The Barnes multiple gamma function}
We assume that $\omega_{1}, \omega_{2}, \cdots, \omega_{n}$ lie on the same
side of some straight line through the origin on the complex plane. 
The Barnes zeta function \cite{bar4} is defined as
  \begin{equation*}
    \zeta_{n}(s,z; {\boldsymbol\omega})
     := \sum_{k_{1},k_{2},\cdots,k_{n}=0}^{\infty}
     \frac{1}{(z+k_{1}\omega_{1}+\cdots+k_{n}\omega_{n})^{s}}.
  \end{equation*}
where ${\boldsymbol \omega}:= (\omega_{1}, \omega_{2},\cdots, \omega_{n})$. 
\par
This is a generalization of the Hurwitz zeta function. 
As a generalization of the formula (\ref{eqn:zgam}), Barnes \cite{bar4}
introduced his multiple gamma functions through
  \begin{equation*}
    \frac{\Gamma_{n}(z, {\boldsymbol\omega})}
      {\rho_{n}({\boldsymbol \omega})}
      := \exp(\zeta'_{n}(0,z;{\boldsymbol\omega})).
  \end{equation*}
where
   \begin{equation*}
     \log\rho_{n}({\boldsymbol \omega})
       :=  -\lim_{z\to 0} \left[\zeta'_{n}(0,z; {\boldsymbol\omega}) 
       + \log z\right].
   \end{equation*}   
It is easy to see that $\Gamma_{n}(z, {\boldsymbol\omega})$ satisfies the
functional relation
  \begin{equation*}
    \frac{\Gamma_{n}(z, {\boldsymbol\omega})}
     {\Gamma_{n}(z+\omega_{i}, {\boldsymbol\omega})}
    =\frac{\rho_{n-1}({\boldsymbol \omega(i)})}
        {\Gamma_{n-1}(z,{\boldsymbol\omega}(i))},
  \end{equation*}
where ${\boldsymbol\omega}(i):= (\omega_{1}, \cdots \omega_{i-1},
\omega_{i+1},\cdots, \omega_{n}).$\par


\subsection{The Vign\'{e}ras multiple gamma function}
As a generalization of the gamma function and the $G$-function, Vign\'{e}ras
\cite{vig} introduced a hierarchy of functions which she called 
``the multiple gamma functions''.

\begin{thm}
There exists a unique hierarchy of functions which satisfy 
 \begin{align*}
      & (1) \quad G_{n}(z+1)=G_{n-1}(z) G_{n}(z), \\
      & (2) \quad G_{n}(1)=1, \\
      & (3) \quad \frac{d^{n+1}}{dz^{n+1}}\log G_{n}(z+1)\geq0
         \quad \mbox{ for } \quad z\geq0,\\
      & (4) \quad G_{0}(z) =z.
  \end{align*}
\label{thm:vig}\end{thm}

Applying Dufresnoy and Pisot's results \cite{duf}, she showed that
these functions satisfying the above properties are uniquely determined
and that $G_{n}(z+1)$ has an infinite product representation
  \begin{align}
     & G_{n}(z+1)=\exp \left[-zE_{n}(1)+\sum_{h=1}^{n-1}
       \frac{p_{h}(z)}{h!}\left(\psi_{n-1}^{(h)}(0)-E_{n}^{(h)}(1)\right)
       \right]\\[8pt]
     & \qquad \times \prod_{{\bold m}\in {\bold N}^{n-1} 
       \times {\bold N}^{*}}
       \left[\left(1+\frac{z}{s({\bold m})}\right)^{(-1)^{n}}
       \exp \left\{ \sum_{l=0}^{n-1}\frac{(-1)^{n-l}}{n-l}
       \left(\frac{z}{s({\bold m})}\right)^{n-l}\right\}\right],
       \nonumber       
  \label{eqn:vig}\end{align}
where
  \begin{align*}
    & E_{n}(z):= \sum_{{\bold m}\in {\bold N}^{n-1} \times {\bold N}^{*}}
      \left[\left\{\sum_{l=0}^{n-1}\frac{(-1)^{n-l}}{n-l}
        \left(\frac{z}{s({\bold m})}\right)^{n-l}\right\}
        +(-1)^{n} \log \left(1+\frac{z}{s({\bold m})}\right)\right],\\
    & \psi_{n-1}(z):=\log G_{n-1}(z+1),\\
    & \frac{e^{tz}-1}{e^{z}-1}=:1+\sum_{k=0}^{\infty}p_{k}(t)\frac{z^{k}}{k!}\\
    & s({\bold m}):=m_{1}+m_{2}+\cdots+m_{n}\quad \mbox{for}
        \quad {\bold m}=(m_{1},m_{2},\cdots, m_{n}),
  \end{align*}
and ${\bold N}^{*}={\bold N}-\{0\}$.\par
 The Vign\'{e}ras multiple gamma function can be regarded as a special
case of the Barnes multiple gamma function. Namely 
  \begin{equation*}
    G_{n}(z)=\Gamma_{n}(z,(1,1\cdots,1))^{(-1)^{n-1}}
      \times(\mbox{the normalization factor}).
  \end{equation*}
In this paper, we will use the word `` the multiple gamma function''
to refer the Vign\'{e}ras multiple gamma function.

\subsection{The $q$-gamma function}
Throughout this paper, we suppose $0<q<1$.
A $q$-analogue of the gamma function was introduced by Jackson \cite{jac1},
\cite{jac2}.
  $$\Gamma (z+1;q)=(1-q)^{-z}
        \prod_{k=1}^{\infty}\left(
            \frac{1-q^{z+k}}{1-q^{k}}
            \right)^{-1}.$$
Askey \cite{ask} pointed out that this function satisfies
a $q$-analogue of the Bohr-Morellup theorem. Namely, $\Gamma(z;q)$
satisfies
    \begin{align*}
        & (1) \quad \Gamma (z+1;q)=[z]_{q} \Gamma(z;q),\\
        & (2) \quad \Gamma (1;q)=1,\\
        & (3) \quad \frac{d^{2}}{dz^{2}}\log\Gamma(z+1;q)\geq0
            \quad \mbox{for}\quad z\geq0,
    \end{align*}
where $[z]_{q}:=(1-q^{z})/(1-q)$. \par
As $q$ tends to $1-0$, $\Gamma(z;q)$ converges $\Gamma(z)$ uniformly with
respect to $z$. A rigorous proof of this fact  was given by Koornwinder
\cite{koo}.\par
 Inspired by Moak's works \cite{moa}, the authors \cite{un} derived a
representation of the $q$-gamma function
  \begin{align}
    \log \Gamma(z;q) & = (z-\frac12) \log \left(
          \frac{1-q^{z}}{1-q} \right)
          + \log q \int_{1}^{z} \xi \frac{q^{\xi}}{1-q^{\xi}}d\xi
          \label{eqn:2}\\
        & + C_{1}(q) + \frac{1}{12} \log q
          + \sum_{k=1}^{m} \frac{B_{2k}}{(2k)!} \left(
          \frac{\log q}{q^{z}-1} \right)^{2k-1} {\widetilde M}_{2k-1}(q^{z})
          \nonumber\\[8pt]
        & + R_{2m}(z;q), \nonumber
  \end{align}
where
  \begin{align*}
    & C_{1}(q) = - \frac{1}{12}\log q -\frac{1}{12}\,\frac{\log q}{q-1}
           + \int _{0}^{\infty} \frac{\overline{B}_{2}(t)}{2}
           \left( \frac{\log q}{q^{t+1}-1} \right)^{2} q^{t+1}dt,\\
    & R_{2m}(z;q) = \int_{0}^{\infty} \frac{\overline{B}_{2m}(t)}{(2m)!}
                \left(
                \frac{\log q}{q^{t+z}-1}
                \right)^{2m} {\widetilde M}_{2m}(q^{t+z})dz,\\[8pt]
    & \overline{B}_{n}(t) := B_{n}(t-[t]) 
       \quad ( \mbox{ $[t]$ denotes the Gauss symbol {\it i.e.} the integral 
         part of } t ) ,
  \end{align*}
and the polynomial ${\widetilde M}_{n}(x)$ is defined by the recurrence 
relation
 \begin{equation}
    {\widetilde M}_{1}(x)=1, \quad
      (x^{2}-x)\frac{d}{dx}{\widetilde M}_{n}(x)
      +nx{\widetilde M}_{n}(x)={\widetilde M}_{n+1}(x).
 \end{equation}
Each term of the formula (\ref{eqn:2}) converges uniformly as $q\to 1-0$.
So we get another proof of the uniformity of the classical limit of
$\log \Gamma(z;q).$

\subsection{The multiple $q$-gamma function}
Recently, one of the authors \cite{nis} constructed the function
$G_{n}(z;q)$ which satisfies a $q$-analogue of the generalized Bohr-Morellup
theorem:

\begin{thm}
There exists a unique hierarchy of functions which satisfy
    \begin{align*}
        & (1) \quad G_{n}(z+1;q)=G_{n-1}(z;q)G_{n}(z;q),\\
        & (2) \quad G_{n}(1;q)=1,\\
        & (3) \quad \frac{d^{n+1}}{dz^{n+1}}\log G_{n+1}(z+1;q)\geq0
            \quad \mbox{for} \quad z\geq0,\\
        & (4) \quad G_{0}(z;q)=[z]_{q}.
    \end{align*}
\label{thm:nis}\end{thm}

We call it ``the multiple $q$-gamma function''. It is given by the following
infinite product representation \cite{nis}
  \begin{align}
    & G_{n}(z+1;q):= (1-q)^{- {z \choose n}} \prod_{k=1}^{\infty}
        \left\{
        \left(
        \frac{1-q^{z+k}}{1-q^{k}}
        \right)^{-k \choose n-1}
        (1-q^{k})^{g_{n}(z,k)}
        \right\} \label{eqn:Gn}
  \end{align}
for $n \geq 1$, where
  $$g_{n}(z,u)={z-u \choose n-1}-{-u \choose n-1}.$$ \par
In the next section, we will derive a representation of the multiple $q$-gamma
function like (\ref{eqn:2}) and consider its classical limit. This limit
formula gives some important properties of the multiple gamma functions.


\section{The Euler-MacLaurin expansion of $G_{n}(z+1;q)$}
By means of the Euler-MacLaurin summation formula
  \begin{align*}
    \sum_{r=M}^{N-1} f(r) 
      &= \int_{M}^{N}f(t)dt + \sum_{k=1}^{n}
        \frac{B_{k}}{k!} \left\{ f^{(k-1)}(N)-f^{(k-1)}(M) 
        \right\}\nonumber\\[8pt] 
      &+ (-1)^{n-1} \int_{M}^{N}\frac{\overline{B}_{n}(t)}{n!}f^{(n)}(t)dt 
        \quad (\mbox{for} \quad f\in C^{n}[M,N]),
  \end{align*}

\vspace{8pt}

we give an expansion formula of the multiple $q$-gamma functions which
we call {\it the Euler-MacLaurin expansion} \cite{un}.  
This formula plays an important role in the following sections.

\begin{prop}
 Suppose $\Re z > -1$ and $m>n$, then
  \begin{align*}
    & \log G_{n}(z+1;q)\\[4pt]
    & \qquad = \left\{ \binom{z+1}{n}+
      \sum_{r=1}^{n} \frac{B_{r}}{r!} \left( - \frac{d}{dz} \right)^{r-1}
      \binom{z}{n-1} \right\} 
      \log \left(\frac{1-q^{z+1}}{1-q}\right)\\[8pt] 
    & \qquad + \sum_{r=1}^{n}\left\{ \left( - \frac{d}{dz} \right)^{r-1}
      \binom{z}{n-1} \right\} \times \int_{1}^{z+1} \frac{\xi^{r}}{r!}
      \frac{q^{\xi} \log q}{1-q^{\xi}} d\xi\\[8pt] 
    & \qquad +  \sum_{j=0}^{n-1}G_{n,j}(z)C_{j}(q) 
      + \sum_{r=1}^{m} \frac{B_{r}}{r!} F_{n,r-1}(z;q) 
      - R_{n,m}(z;q),
  \end{align*}
where {\allowdisplaybreaks
  \begin{align*}
    & F_{n,r-1}(z;q) := \left[\frac{d^{r-1}}{dt^{r-1}}
      \left\{\binom{-t}{n-1}\log \left( \frac{1-q^{z+t}}{1-q^{z+1}}
      \right)
     \right\}\right]_{t=1},\\[8pt]
    &  C_{j}(q):= -\sum_{r=1}^{n+1}\frac{B_{r}}{r!}f_{j+1,r-1}(q)\\
    &  \qquad + \frac{(-1)^{n}}{(n+1)!}\int_{1}^{\infty}
           \left[ \overline{B}_{n+1}(t)
           \left\{\frac{d^{n+1}}{dt^{n+1}}
              \left\{t^{j}\log
                \left(\frac{1-q^{t}}{1-q} \right)
              \right\}
           \right\} \right]dt,\\[8pt]
    &  f_{j+1,r-1}(q):=
         \left[\frac{d^{r-1}}{dt^{r-1}}
            \left\{t^{j} \log
            \left( \frac{1-q^{t}}{1-q} \right)
            \right\} \right]_{t=1},\\[8pt]
    &   R_{n,m}(z;q):=  \frac{(-1)^{m-1}}{m!}
         \int_{1}^{\infty}
         \left[ \overline{B}_{m}(t)
         \left\{\frac{d^{m}}{dt^{m}}
            \left\{\binom{-t}{n-1} \log
            \left(\frac{1-q^{z+t}}{1-q^{z+1}}\right)
            \right\}
         \right\} \right]dt.
  \end{align*}}
The polynomial $G_{n,j}(z)$ is introduced through
  $$\binom{z-u}{n-1}= \sum_{j=0}^{n-1}G_{n,j}(z)u^{j}.$$
\label{prop;1}\end{prop}

\begin{pf}
From (2.7) and the definition of $G_{n,j}(z)$, we obtain
  \begin{align}
    \log G_{n}(z+1;q) 
      &= -\binom{z}{n}\log(1-q) - \sum_{k=1}^{\infty}
        \binom{-k}{n-1}\log (1-q^{z+k})\label{eqn;1}\\[8pt] 
      &+ \sum_{j=0}^{n-1} G_{n,j}(z)\{\sum_{k=1}^{\infty}
        k^{j}\log(1-q^{k})\}\nonumber.
  \end{align}
\vspace{8pt}
Applying the Euler-MacLaurin summation formula, we have
{\allowdisplaybreaks
  \begin{align}
    & \sum_{k=1}^{\infty}\binom{-k}{n-1}
      \log(1-q^{z+k})\label{eqn;2}\\[8pt]
    & \qquad =\int_{1}^{\infty} 
      \binom{-t}{n-1}\log(1-q^{z+t})dt\nonumber\\[8pt] 
    & \qquad - \sum_{r=1}^{m}\frac{B_{r}}{r!} \left\{
      \left(\frac{d}{dt}\right)^{r-1} \binom{-t}{n-1}
      \right\}\Bigg\vert_{t=1}\log(1-q^{z+1})\nonumber\\[8pt] 
    & \qquad -  \sum_{r=1}^{m}\frac{B_{r}}{r!}F_{n,r-1}(z;q) 
      + R_{n,m}(z;q).\nonumber
  \end{align}}
Let $L_{r}(z)$ and $L_{r}$ be
  $$L_{r}(z):= \frac{Li_{r}(q^{z})}{\log^{r-1}q},
    \quad L_{1}(z):= -\log(1-q^{z}), \quad
  L_{r}:=L_{r}(q),$$
where $Li_{r}(z)$ is Euler's polylogarithm
  $$Li_{r}(z):=\sum_{k=1}^{\infty}\frac{z^{k}}{k^{r}}.$$ 
From
  $$\log(1-q^{z+t})=-\sum_{l=1}^{\infty}\frac{q^{(t+z)l}}{l},$$
it follows that
  \begin{align}
    & \int_{1}^{\infty}\binom{-t}{n-1}\log(1-q^{t+z})dt \label{eqn;3}\\[8pt] 
    & \qquad \quad = \sum_{j=0}^{n-1}
    \frac{(-1)^{j}\sideset{_{n-1}}{_j}S}{(n-1)!}
    \sum_{r=0}^{j}\frac{(-1)^{r}j!}{(j-r)!}L_{r+2}(z+1),\nonumber
  \end{align}

where $\sideset{_{n}}{_j}S$ is the Stirling number of the first kind
defined by
 $$ \sum_{j=0}^{n}\sideset{_{n}}{_j}S u^{j}
    = [u]_{n},$$
where $[u]_{n}=u(u-1)\cdots(u-n+1).$
Substituting (\ref{eqn;3}) to (\ref{eqn;2}), we have 
{\allowdisplaybreaks
  \begin{align}
    &  \sum_{k=1}^{\infty}\binom{-k}{n-1}\log(1-q^{z+k})\label{eqn;4}\\[4pt]
    &  \quad = \sum_{j=0}^{n-1}
      \frac{(-1)^{j}\sideset{_{n-1}}{_j}S}{(n-1)!}
      \sum_{r=0}^{j} \frac{(-1)^{r}j!}{(j-r)!}L_{r+2}(z+1)\nonumber\\[8pt]
    &  \quad    + \sum_{r=1}^{m} \frac{B_{r}}{r!}
      \left\{ \left(
        \frac{d}{dt} \right)^{r-1} \binom{-t}{n-1}
      \right\}\Bigg\vert_{t=1}L_{1}\nonumber\\[8pt]
    &  \quad - \sum_{r=1}^{m} \frac{B_{r}}{r!}F_{n, r-1}(z;q)
       + R_{n,m}(z;q).\nonumber
  \end{align}} 

Similarly, 
  \begin{align}
    &   \sum_{k=1}^{\infty} k^{j}\log(1-q^{k})\label{eqn;5}\\
    &  \qquad    = \sum_{r=0}^{j}
      \frac{(-1)^{r}j!}{(j-r)!} L_{r+2}
      + \sum_{r=1}^{j}\frac{B_{r}}{r!}
        \frac{j!}{(j+1-r)!}L_{1} + C_{j}(q).\nonumber
  \end{align}

$L_{r}(z)$ and $L_{r}$ cause divergence as $q \to 1-0$, but we can 
prove that these divergent terms vanish. In order to simplify 
such terms, let us express $L_{r}(z)$ as the sum of $L_{r}$ and 
convergent terms. 

\begin{lem}
  \begin{eqnarray*}
    & & L_{l+1}(z) = \frac{z^{l}}{l!}\log
      \left(\frac{1-q^{z}}{1-q}\right)\\
    & &\qquad  + \sum_{r=0}^{l}\frac{(z-1)^{l-r}}{(l-r)!}L_{r+1}
      +  \sum_{r=1}^{l}\frac{(-1)^{l}z^{l-r}}{(l-r)!}T_{r}(z),
  \end{eqnarray*}
where
  $$T_{r}(z)=\int_{1}^{z}\frac{\xi^{r}}{r!}
    \frac{q^{\xi}\log q}{1-q^{\xi}} d\xi. $$
\label{lem;2}\end{lem}    

\begin{pf}
By partial integration, we have 
  \begin{equation}
    \int_{1}^{z}\xi^{n}\frac{q^{\xi}\log q}{1-q^{\xi}} d\xi
     = \sum_{r=0}^{n}\frac{n!(-1)^{r}}{(n-r)!}\frac1{\log^{r}q}
     \left\{z^{n-r}Li_{r+1}(q^{z})-Li_{r+1}(q)\right\}.
  \label{eqn;5.5}\end{equation}

Let ${\cal L}_{k}(z)$ be
  \begin{equation*}
    {\cal L}_{k+1}(z)=\frac{Li_{k+1}(q^{z})-Li_{k+1}(q)}{\log^{k}q}.
  \end{equation*}

Then (\ref{eqn;5.5}) reads
  \begin{align}
    {\cal L}_{n+1}(z)&= -\sum_{k=0}^{n-1}
    \frac{(-1)^{n-k}z^{n-k}}{(n-k)!}{\cal L}_{k+1}(z)\\[8pt]
    & - \sum_{k=0}^{n-1}\frac{(-1)^{n-k}(z^{n-k}-1)}{(n-k)!}L_{k+1}
    + (-1)^{n}T_{n}(z).\nonumber
  \label{eqn;6}\end{align}

\vspace{8pt}

In order to write ${\cal L}_{n+1}(z)$ with $L_{n+1}$ and $T_{k}(z)$, we use
induction on n. It can be seen that
  \begin{equation*}
    {\cal L}_{n+1}(z) = \frac{z^{n}{\cal L}_{1}(z)}{n!}
      +\sum_{k=0}^{n-1}\frac{(z-1)^{n-k}}{(n-k)!}L_{k+1}
      +\sum_{k=1}^{n}\frac{(-1)^{k}z^{n-k}}{(n-k)!}T_{k}(z),
   \end{equation*}
which is equivalent to the claim of Lemma \ref{lem;2}.
\end{pf}

Applying (\ref{eqn;4}), (\ref{eqn;5}) and Lemma \ref{lem;2} to (\ref{eqn;1}),
we obtain

{\allowdisplaybreaks
  \begin{align}
    & \log G_{n}(z+1;q) \label{eqn;8}\\[4pt]
    &  \quad = - \sum_{j=0}^{n-1}\frac{(-1)^{j} \sideset{_{n-1}}{_j}S}{(n-1)!}
       \sum_{l=1}^{j}\frac{(-1)^{l}j!}{(j-l)!}
       \sum_{r=0}^{l}\frac{z^{l-r}}{(l-r)!}L_{r+2}\nonumber\\[8pt]
    &  \qquad + \sum_{j=0}^{n-1}G_{n,j}(z) 
    \sum_{r=0}^{j}\frac{(-1)^{r}j!}{(j-r)!}
       L_{r+2}\nonumber\\[8pt]
    &  \qquad +  \left\{\binom{z}{n} - \sum_{r=1}^{m}\frac{B_{r}}{r!}
       \left(\frac{d}{dt} \right)^{r-1}\binom{-t}{n-1}\Bigg\vert_{t=1}
       \right.\nonumber\\[8pt]
    &  \qquad \qquad - \sum_{j=0}^{n-1}
      \frac{(-1)^{j}\sideset{_{n-1}}{_j}S}{(n-1)!}
      \sum_{l=0}^{j}\frac{(-1)^{l}j!}{(j-l)!}
      \frac{z^{l+1}}{(l+1)!}\nonumber\\[8pt]
    &  \left.\qquad \qquad \qquad \qquad
      + \sum_{j=0}^{n-1}G_{n,j}(z)\sum_{r=0}^{j}
      \frac{B_{r}}{r!} \frac{j!}{(j+1-r)!}
      \right\} L_{1}\nonumber\\[8pt]
    & \qquad +\left\{-\sum_{j=0}^{n-1}
      \frac{(-1)^{l}\sideset{_{n-1}}{_j}S}{(n-1)!}
      \sum_{l=0}^{j}\frac{(-1)^{l}j!}{(j-l)!}
      \frac{(z+1)^{l+1}}{(l+1)}\right.\nonumber\\[8pt]
    & \left. \qquad \qquad + \sum_{j=0}^{n-1} \frac{B_{j+1}}{(j+1)!}
      \left(\frac{d}{dt}\right)^{j}
      \binom{-t}{n-1}\Bigg\vert_{t=1}
      \right\}
      \times \log \left(
        \frac{1-q^{z+1}}{1-q}
        \right)\nonumber\\[8pt]
    & \qquad -\sum_{j=0}^{n-1}\frac{(-1)^{j}\sideset{_{n-1}}{_j}S}{(n-1)!}
      \sum_{l=0}^{j}\frac{(-1)^{l}j!}{(j-l)!}
      \sum_{r=0}^{l}\frac{(-1)^{r}(z+1)^{l-r}}{(l-r)!}
      T_{r+1}(z+1)\nonumber\\[8pt]
    &  \qquad + \sum_{r=0}^{m}\frac{B_r}{r!}F_{n.r-1}(z;q)
      + \sum_{j=1}^{n-1}G_{n,j}(z)C_{j}(q)\nonumber\\[4pt]
    &  \qquad -R_{n,m}(z;q).\nonumber
  \end{align}}

We prove that the coefficients of the divergent terms in (\ref{eqn;8})
vanish. By the following lemma, we can see that    
the first term and the second term in (\ref{eqn;8}) are canceled out.

\begin{lem}
\begin{align*}
 & \sum_{j=0}^{n-1}\frac{(-1)^{j} \sideset{_{n-1}}{_j}S}{(n-1)!}
     \sum_{l=1}^{j}\frac{(-1)^{l}j!}{(j-l)!}
     \sum_{r=0}^{l}\frac{z^{l-r}}{(l-r)!}L_{r+2}\\[8pt]
 &  \qquad \qquad = \sum_{j=0}^{n-1}G_{n,j}(z)
     \sum_{r =0}^{j}\frac{(-1)^{r}j!}{(j-r)!}
     L_{r+2}.
\end{align*}
\end{lem}
\vspace{8pt}

\begin{pf}
We consider the exponential generating functions of the
coefficients of $L_{r+2}$ in both sides.
The generating function of the left hand side is
{\allowdisplaybreaks
  \begin{align*}
    & \sum_{j=0}^{n-1}\frac{(-1)^{j}\sideset{_{n-1}}{_j}S}{(n-1)!}
     \sum_{l=0}^{j}\frac{(-1)^{l}j!}{(j-l)!}
    \sum_{r=0}^{l}\frac{z^{l-r}}{(l-r)!}\frac{u^{r}}{r!}\\[8pt]
    & \quad = \sum_{j=0}^{n-1}\frac{(-1)^{j}\sideset{_{n-1}}{_j}S}{(n-1)!}
        (1-z-u)^{j}\\[8pt]
    & \quad = \binom{z+u-1}{n-1}.
 \end{align*}}
\vspace{8pt}

On the other hand, the generating function of the right hand side is
{\allowdisplaybreaks
  \begin{align*}
   & \sum_{j=0}^{n-1} G_{n,j}(z)
     \sum_{r=0}^{j}\frac{(-1)^{r}j!}{(j-r)!}
     \frac{u^{r}}{r!}\\[8pt]
   & \quad  = \sum_{j=1}^{n-1}G_{n,j}(z)(1-u)^{j}\\[8pt]
   & \quad = \binom{z+u-1}{n-1}.
 \end{align*}}
\vspace{8pt}

Therefore, the coefficients of $L_{r+2}$ in both sides coincide.
Thus, the claim follows.
\end{pf}  

Next, we prove the coefficient of $L_{1}$ vanishes. 
Since
  \begin{equation*}
    \sum_{j=0}^{n-1}\frac{(-1)^{j}\sideset{_{n-1}}{_j}S}{(n-1)!}
     \sum_{l=0}^{j}\frac{(-1)^{l}j!}{(j-l)!} \frac{z^{l+1}}{(l+1)!}
    = \int_{0}^{z}\binom{t-1}{n-1}dt
  \end{equation*}
and
{\allowdisplaybreaks
  \begin{align*}
    & \sum_{j=0}^{n-1}G_{n,j}(z)\sum_{r=0}^{j}
      \frac{B_{r}}{r!}\frac{j!}{(j+1-r)!}\\[8pt]
    & \qquad   =  \sum_{r=0}^{n-1}\frac{B_{r}}{r!}
      \sum_{j=r}^{n-1}G_{n,j}(z)
      \frac{j!}{(j+1-r)!}\\[8pt]
    & \qquad = \sum_{r=0}^{n-1}\frac{B_{r}}{r!}
       \sum_{j=r}^{n-1}G_{n,j}(z)
       \left\{\left(\frac{d}{dt}\right)^{r-1} t^j\right\}
       \Bigg\vert_{t=1}\\[8pt]
    & \qquad   =  \sum_{r=0}^{n-1}\frac{B_{r}}{r!}
      \left\{\left(\frac{d}{dt}\right)^{r-1}
      \binom{z-t}{n-1} \right\}\Bigg\vert_{t=1}\\[8pt]
    & \qquad = \sum_{r=0}^{n-1}\frac{B_{r}}{r!}
      \left\{\left(-\frac{d}{dz}\right)^{r-1}
      \binom{z-1}{n-1}\right\},
  \end{align*}}

we have
  \begin{align*}
   & (\mbox{the coefficient of $L_{1}$})\\[4pt]
   &  \quad = \binom{z}{n} - \int_{0}^{z} \binom{t-1}{n-1}dt
     + \sum_{r=0}^{m}\frac{B_{r}}{r!}
     \left( - \frac{d}{dz}\right)^{r-1}\binom{z-1}{n-1}\\[8pt]
   & \qquad  -\sum_{r=0}^{m}\frac{B_{r}}{r!}
     \left( - \frac{d}{dz}\right)^{r-1}\binom{z-1}{n-1}
     \Bigg\vert_{z=0}.
  \end{align*}
 
We prove that the right hand side of the above formula is
equal to zero. In a formal sense,
  \begin{equation*}
    \left( e^{-\frac{d}{dz}}-1\right)^{-1}
      +\left(\frac{d}{dz}\right)^{-1}
    = \sum_{r=1}^{\infty}\frac{B_{r}}{r!}
      \left(-\frac{d}{dz}\right)^{r-1}.
  \end{equation*}

Imposing the boundary condition at $z=0$, we make 
the both sides above act on $\binom{z-1}{n-1}$. Then,
  \begin{equation*}
    \left(e^{-\frac{d}{dz}}-1\right)^{-1}\binom{z-1}{n-1}
       + \int_{0}^{z} \binom{t-1}{n-1} dt
       + \sum_{r=1}^{n}\frac{B_{r}}{r!}
       \left[\left(-\frac{d}{dt}\right)^{r-1}
       \binom{t-1}{n-1}
    \right]_{t=0}^{t=z} =0 
  \end{equation*}
because $\binom{t-1}{n-1}$ is a polynomial of $(n-1)$-degree.
Since $F(z):=-\binom{z}{n}$ satisfies
  \begin{equation*}
    F(0)=0 \qquad \binom{z-1}{n-1}=F(z-1)-F(z),
  \end{equation*}
it can be seen that
  \begin{equation*}
    \left(e^{-\frac{d}{dz}}-1\right)^{-1}\binom{z-1}{n-1}
    =-\binom{z}{n}.
  \end{equation*}

Thus we have the formula
  \begin{equation}
    -\binom{z}{n}
     + \int_{0}^{z} \binom{t-1}{n-1} dt
     + \sum_{r=1}^{n}\frac{B_{r}}{r!}
   \left[\left(-\frac{d}{dt}\right)^{r-1}
   \binom{t-1}{n-1}
   \right]_{t=0}^{t=z} =0, \label{eqn;9}
  \end{equation}
which shows that the coefficient of $L_{1}$ is equal to zero.
Hence we have proved that the all coefficients of $L_{r}$
vanish in (\ref{eqn;8}).\par
Finally, we calculate the coefficients of 
$\log\left(\frac{1-q^{z+1}}{1-q}\right)$ and $T_{r}(z+1)$. 
Using the formula (\ref{eqn;9}), we have
\begin{align}
 & \left(\mbox{the coefficient of 
     $\log\left(\frac{1-q^{z+1}}{1-q}\right)$ in (\ref{eqn;1})}
     \right)\label{eqn;10}\\[4pt]
 & \qquad \qquad  \qquad= \binom{z+1}{n}+ \sum_{r=1}^{n}
     \frac{B_{r}}{r!}\left(
     - \frac{d}{dt} \right)^{r-1}
     \binom{z}{n-1}.\nonumber 
\end{align}

In order to calculate the coefficients of $T_{r}(z+1)$
in (\ref{eqn;8}), we note that 
  \begin{align*}
    \binom{z-u}{n-1}
     &= \sum_{j=0}^{n-1}\frac{\sideset{_{n-1}}{_j}S}{(n-1)!}
       (z-u)^{j}\\[8pt]
     &= \sum_{j=0}^{n-1}\frac{(-1)^{j}\sideset{_{n-1}}{_j}S}{(n-1)!}
       \sum_{l=0}^{j}\frac{(-1)^{l}j!}{(j-l)!}
       \sum_{r=0}^{l}\frac{(-1)^{r}(z+1)^{l-r}}{(l-r)!}
       \frac{u^{r}}{r!}.
  \end{align*}

Using this identity, we have
  \begin{align}
     & \left(\mbox{the coefficient of $T_{r}(z+1)$ in (\ref{eqn;8})}
      \right)\label{eqn;11}\\[4pt]
     &   \qquad \qquad= \left(\frac{d}{dt}\right)^{r-1}
        \binom{z-u}{n-1}\Bigg\vert_{u=0}
        = \left(-\frac{d}{dz}\right)^{r-1}
        \binom{z}{n-1}.\nonumber
  \end{align}

Substituting (\ref{eqn;10}) and (\ref{eqn;11}) in (\ref{eqn;8}),
we obtain Proposition 1.1.
\end{pf}


\section{The classical limit of $G_{n}(z+1;q)$ and 
the asymptotic expansion of $G_{n}(z+1)$}
 
In this section, we study the classical limit of $G_{n}(z+1;q)$
using the Euler-MacLaurin expansion.  We will see that this limit
formula gives an asymptotic expansion for the multiple gamma functions, 
which is a generalization of the Stirling formula for the gamma function.


\subsection{The classical limit of $G_{n}(z+1;q)$}
First we consiser the classical limit of the Euler-MacLaurin expansion
in the domain  $\{z \in {\bold C}|\Re z >-1 \}$.

\begin{prop}
Suppose $\Re z > -1$ and $m>n$.
  \begin{align*}
   & \lim_{q\to 1-0} \log G_{n}(z+1;q)\\[4pt]
     & \quad = \left\{ \binom{z+1}{n}+ \sum_{r=1}^{n}
         \frac{B_{r}}{r!}
         \left( - \frac{d}{dz} \right)^{r-1}
         \binom{z}{n-1} \right\}
          \log(z+1) \\[8pt]
    & \quad - \sum_{r=1}^{n}\left\{
          \left( - \frac{d}{dz} \right)^{r-1}
          \binom{z}{n-1}\right\}\frac{1}{r!r}
          \left\{(z+1)^{r}-1\right\}\\[8pt]
    & \quad - \sum_{j=0}^{n-1}G_{n,j}(z)C_{j}
         + \sum_{r=1}^{m} \frac{B_{r}}{r!} F_{n,r-1}(z)\\[8pt]
    & \quad - R_{n,m}(z), 
    \end{align*}
\vspace{8pt}

where
{\allowdisplaybreaks
 \begin{align*}
   & C_{j}:= -\sum_{r=1}^{n+1}\frac{B_{r}}{r!} 
       \left(\frac{d}{dt}\right)^{r-1}
       \{t^{j} \log t \}\Bigg\vert_{t=1}\\
   & \qquad \qquad \qquad \qquad
       + \frac{(-1)^{n}}{(n+1)!}\int_{1}^{\infty}
       \overline{B}_{n+1}(t)\left(
       \frac{d}{dt}\right)^{n+1}
       \left\{t^{j}\log t\right\} dt\\[8pt]
   &  F_{n,r-1}:=\left(\frac{d}{dt}\right)^{r-1}
      \left\{\binom{-t}{n-1}
        \log\left(\frac{z+t}{z+1}\right)
        \right\}\Bigg\vert_{t=1},\\[8pt]
   & R_{n,m}(z):=\frac{(-1)^{m-1}}{m!}\int_{1}^{\infty}
       \overline{B}_{m}(t) \left(
       \frac{d}{dt}\right)^{m}
       \left\{\binom{-t}{n-1}
       \log \left(\frac{z+t}{z+1}\right)\right\}dt.
  \end{align*}}
\vspace{8pt}

Furthermore this convergence is uniform on any compact
set in $\{z \in {\bold C}|\Re z > -1 \}$.
\label{prop;21}\end{prop}

\begin{pf}
Taking Proposition \ref{prop;1} into acount, we must show that 
{\allowdisplaybreaks
  \begin{align}
     & \lim_{q\to1-0} \log \left(\frac{1-q^{z+1}}{1-q}\right) =
       \log \left(z+1\right),\\[8pt]
     & \lim_{q\to 1-0} \int_{1}^{z+1}\frac{\xi^{r}}{r!}
         \frac{q^{\xi}\log q}{1-q^{\xi}}d\xi =
         -\int_{1}^{\infty}\frac{\xi^{r-1}}{r!}d\xi
         =- \frac{1}{r!r} \left\{\left( z+1 \right)^{r}-1\right\},\\[8pt]
     & \lim_{q\to 1-0} F_{n,r-1}(z;q) = F_{n,r-1}(z),\\[4pt]
     & \lim_{q\to 1-0} R_{n,m}(z;q) = R_{n,m}(z)\label{eqn:(4)},\\[4pt]
     &  \lim_{q\to 1-0} C_{j}(q) = C_{j},
  \end{align}}
and further have to show that this convergence is uniform. Here we prove 
only (\ref{eqn:(4)}). The other formulas can be verified in a similar way.\par
Since
  \begin{equation*}
    \lim_{q\to1-0}  \left(\frac{d}{dt}\right)^{r-1}\left\{
      \binom{-t}{n-1}\log\left(\frac{1-q^{z+t}}{1-q^{z+1}}\right)
      \right\}
    = \left(\frac{d}{dt}\right)^{r-1}
      \left\{ \binom{-t}{n-1} \log \left(
      \frac{z+t}{z+1}\right)\right\},
  \end{equation*}
in order to show (\ref{eqn:(4)}), it is sufficient to prove that the 
procedure of taking the classical limit commutes with the integration. 
 Let us introduce polynomials $M_{r}(x)$ through
  \begin{equation*}
    \frac{d^{r}}{dz^{r}}\log (1-q^{z+t})
        = - \left( \frac{\log q}{1-q^{z+t}}\right)^{r}
        q^{z+t} M_{r}(q^{z+t})
 \end{equation*}
(cf.[18]). They satisfy the recurrence relation.
  \begin{equation*}
      M_{1}(x)=1, \quad (x^{2}-x) \frac{d}{dx}M_{n}(x)+\{(r-1)x+1\}M_{r}(x)
      =M_{r+1}(x)
  \end{equation*}
and $M_{r}(1)=(r-1)!$. Using this we have 
{\allowdisplaybreaks
  \begin{align*}
    & \int_{1}^{\infty}\overline{B}_{m}(t)
      \left(\frac{d}{dt}\right)^{m}
      \binom{-t}{n-1} \log
      \left(\frac{1-q^{z+t}}{1-q^{z+1}}\right)dt\\[8pt]
    & \qquad =\frac{1}{(n-1)!}\sum_{j=1}^{n-1}
      \left\{(-1)^{j}\sideset{_{n-1}}{_j}S \sum_{l=0}^{j}
      \binom{m}{l}[j]_{l}\right.\\[8pt]
    & \left.\qquad\qquad \times \int_{1}^{\infty}\overline{B}_{m}(t)t^{j-l}
      \left(\frac{\log q}{1-q^{t+z}}\right)^{m-l}q^{t+z}
      M_{m-l}(q^{t+z})dt\right\}.
  \end{align*}}
Therefore we have to show 
  \begin{align}
    & \lim_{q\to1-0}\int_{1}^{\infty}\overline{B}_{m}(t)
       \left\{t^{j-l}\left(\frac{\log q}{1-q^{t+z}}\right)^{m-l}
       q^{t+z}M_{m-l}(q^{t+z})\right\}dt
       \label{eqn:2*}\\[8pt]
    &  \qquad = \int_{1}^{\infty} \lim_{q\to1-0}\overline{B}_{m}(t)
       \left\{t^{j-l}\left(\frac{\log q}{1-q^{t+z}}\right)^{m-l}
       q^{t+z}M_{m-l}(q^{t+z})\right\}dt.
       \nonumber
  \end{align}
\vspace{8pt}

\begin{lem}
Suppose  that $\alpha \in {\bold N}$ be fixed
and that $y_0$ and $y_1$ are fixed constants such that
$-1 < y_0 < y_1$. Then there exists a constant $C$
depending on $y_0$ and $y_1$ such that
  \begin{equation*}
    0 < \left(\frac{\log q}{q^{t+y}-1} \right)^{\alpha}
      < \frac{C}{(t+y)^{\alpha}}
    \quad \mbox{for} \quad y \in [y_0, y_1], \quad
      1<t<\infty, \quad 0<q<1.
  \end{equation*}
\label{lem:22}\end{lem}

\begin{pf}
Letting $q=e^{-\delta}$, we define 
  \begin{equation*}
    \varphi(t,y,\delta):=
      \begin{cases}
        \left(\frac{(t+y)\delta}{e^{-\delta(t+y)}-1}\right)^{\alpha}
         e^{-\delta (t+y)}
         & \qquad \delta > 0\\
        \qquad 1 & \qquad \delta = 0.
      \end{cases}
  \end{equation*}
Then, $\varphi(t,y,\delta)$ is positive, continuous and bounded
on
  \begin{equation*}
    D:=\left\{(t,y,\delta)\Bigg\vert 1 \leq t \leq +\infty,\quad
        \delta \geq 0,\quad y_0 \leq y \leq y_1
        \right\}.
  \end{equation*}
So we can put $C:= \underset{D}{\sup}\,\varphi(t,y,\delta).$
\end{pf}

Suppose $-1<y_{0} \leq \Re z \leq y_{1}.$ By Lemma \ref{lem:22}, there 
exists a constant $C$ depending only on $y_{0}$ and $y_{1}$ such that
  \begin{equation*}
    \Bigg\vert t^{j-l} \left(
      \frac{\log q}{1-q^{t+z}}\right)^{m-l}
      q^{t+z} M_{m-l}(q^{t+z})
    \Bigg\vert
    \leq t^{j-l}\frac{C}{(t+\Re z)^{m-l}}
  \end{equation*}
for $1\leq t <+\infty$, $0 \leq q <1.$\par
Since $m-j \geq2,$ the right hand side above is integrable over
$1\leq t < +\infty$.  Therefore, Lebesgue's convergence theorem
ensures (\ref{eqn:2*}). We have thus proved the limit formula
(\ref{eqn:(4)}). Next we show that this convergence is uniform on any 
compact set in the domain $\{ z \in {\bold C}|\Re z >-1\}$.\par
Put 
  \begin{equation*}
    \Phi(z,q):=\int_{1}^{\infty}
      \overline{B}_{m}(t)t^{j-l}\left(
      \frac{\log q}{1-q^{z+t}}
      \right)^{m-l}q^{t+z}M_{m-l}(q^{z+t})dt.
  \end{equation*}
From the consideration above, we see that there exists a constant 
$C$ depending only on $y_{0}$ and $y_{1}$ such that
  \begin{equation*}
    |\Phi(z,q)| \leq
      C \int_{1}^{\infty}\overline{B}_{m}(t)\frac{dt}{t^{2}}
      \quad \mbox{for} \quad \Re z \in [y_0,y_1]
      \quad \mbox{and} \quad 0 < q \leq 1.
  \end{equation*}
Hence $\{\Phi(z,q)|0 < q \leq 1\}$ is a uniformly bounded family
of functions and
  \begin{equation*}
    \Phi(z,q)\to\int_{1}^{\infty}\overline{B}_{m}(t)
      \frac{t^{j-l}(-1)^{m-l}(m-l-1)!}{(t+z)^{m-l}}dt
    \quad \mbox{as} \quad q \to 1-0.
  \end{equation*}  
By Vitali's convergence theorem, this convergence is uniform
on any compact set in the domain $\{ z \in {\bold C}|\Re z >-1\}$.
\end{pf}

The constant $C_{j}$ in Proposition 4.1 can be expressed in terms
of the special value of the Riemann zeta function.

\begin{lem}
  \begin{equation*}
    C_{j}=-\zeta'(-j)-\frac{1}{(j+1)^{2}}.
  \end{equation*}
\label{lem;23}\end{lem}

\begin{pf}
From the definition of $\zeta(s)$, 
 \begin{equation*}
   \zeta ' (s) = - \sum_{k=1}^{\infty}
     \frac{\log k}{k^{s}}\quad
     \mbox{for}\quad \Re s > 1.
 \end{equation*}    
By the Euler-MacLaurin summation formula, we obtain
{\allowdisplaybreaks
  \begin{align}
    \zeta'(s) &= - \frac{1}{(s-1)^{2}}
        + \sum_{r=1}^{n} \frac{B_{r}}{r!}
        \left(\frac{d}{dt}\right)^{r-1}
        \left\{t^{-s}\log t\right\}\Bigg\vert_{t=1}\label{eqn;23}\\[8pt]
        & \qquad -\frac{1}{n!}\int_{1}^{\infty}\overline{B}_{n}(t)
        \left(\frac{d}{dt}\right)^{n}
        \left\{t^{-s}\log t\right\}dt.\nonumber
  \end{align}}
Since
  \begin{equation*} 
    \left(\frac{d}{dt}\right)^{n} \{t^{-s}\log t\}
     = -\frac{\partial}{\partial s}[-s]_{n}t^{-s-n}
       +[-s]_{n}t^{-s-n}\log t,
  \end{equation*}
(\ref{eqn;23}) can be analytically continued to $\{z \in {\bold C}|
\Re z > -n+1\}$. So if we put $s=-j$, $n=j+2$, then we obtain
{\allowdisplaybreaks
  \begin{align*}
    \zeta'(-j) &= -\frac{1}{(j+1)^{2}}+
      \sum_{r=1}^{j+2}\frac{B_r}{r!}
      \left(\frac{d}{dt}\right)^{r-l}
      \left\{t^{j} \log t\right\}
      \Bigg\vert_{t=1}\\[8pt]
      &\qquad - \frac{(-1)^{j+1}}{(j+2)^{2}} \int_{1}^{\infty} 
      \overline{B}_{j+2}(t)\left(\frac{d}{dt}\right)^{j+2} 
      \{t^{j} \log t \}dt \\[8pt]
      &=  -\frac{1}{(j+1)^{2}}-C_{j}.
  \end{align*}}
\end{pf}

Next we prove that the limit function in Proposition 4.1 coincides with
the multiple gamma funciton.

\begin{thm}
Suppose $m>n$. Then, as $q\to1-0$, $G_{n}(z+1;q)$ converges to $G_{n}(z+1)$ 
uniformly on any compact set in the domain ${\bold C}\backslash 
{\bold Z}_{<0}$ and
{\allowdisplaybreaks
  \begin{align} 
   &\log G_{n}(z+1)\label{eqn;22}\\[4pt]
       & \quad=  \left\{ \binom{z+1}{n}+ \sum_{r=1}^{n}
         \frac{B_{r}}{r!}
         \left( - \frac{d}{dz} \right)^{r-1}
         \binom{z}{n-1} \right\}
          \log(z+1) \nonumber\\[8pt]
       & \quad - \sum_{r=1}^{n} \left\{
          \left( - \frac{d}{dz} \right)^{r-1}
          \binom{z}{n-1}\right\}
          \times \frac{1}{r!r}
          \left\{(z+1)^{r}-1\right\}\nonumber\\[8pt]
       & \quad-  \sum_{j=0}^{n-1}G_{n,j}(z)
          \left\{\zeta'(-j)+\frac{1}{(j+1)^{2}}\right\}
         + \sum_{r=1}^{m} \frac{B_{r}}{r!} F_{n,r-1}(z)\nonumber\\[8pt]
       & \quad - R_{n,m}(z).\nonumber
    \end{align}} 
\end{thm}

\begin{pf}
In the domain $\{z \in {\bold C}|\Re z > - 1\}$, we have already proved 
the existence of the limit function and uniformity of the convergence. 
Let us put
  \begin{equation*}
    \tilde{G}_{n}(z+1):=\lim_{q \to 1-0}G_{n}(z+1;q).
  \end{equation*}  
Because of the uniformity of the convergence, we have particularly
  \begin{equation*}
    \lim_{q \to 1-0} \left[ \left( \frac{d}{dz}\right)^{n+1} 
      \left\{\log G_{n}(z+1;q)\right\} \right]
    = \left(\frac{d}{dz}\right)^{n+1} 
      \left\{\log \tilde{G}_{n}(z+1) \right\}
  \end{equation*}    
so that, from Theorem 2.2, $\tilde{G}_{n}(z+1)$ satisfies the conditions in
Theorem 2.1. Namely
  \begin{align*}
    & \mbox{(1)} \quad \tilde{G}_{n}(z+1) 
      =\tilde{G}_{n-1}(z)\tilde{G}_{n}(z)\\[4pt]
    & \mbox{(2)} \quad \left(\frac{d}{dz}\right)^{n+1} 
      \log \tilde{G}_{n}(z+1) \geq 0 
      \quad \mbox{for} \quad z \geq 0\\[4pt] 
    & \mbox{(3)}\quad\tilde{G}_{n}(1)=1\\[4pt] 
    & \mbox{(4)}\quad \tilde{G}_{0}(z+1)=z+1.
  \end{align*}
Since a hierarchy of such functions is uniquely determined, so  
$G_{n+1}(z+1)=\tilde{G}_{n}(z+1)$ in $\{z\in {\bold C}|\Re z >-1\}$. 
Thus the claim of the theorem in the case that $\Re z>-1$ has been 
proved for $\{z\in {\bold C}|\Re z >-1\}$.\par 
Next, we show that in $\{z\in {\bold C}|\Re z \leq -1,z\ne -1\}$,
  \begin{equation*}
    G_{n}(z+1;q) \to G_{n}(z+1) \quad \mbox{as} \quad q\to1-0
  \end{equation*}
and that the convergence is uniform on any compact set in this domain.
For the proof, we use induction on $n$.\par 
The case that $n=1$ was considered by Koornwinder \cite{koo}. 
Let $K$ be a compact set in $\{z \in {\bold C}|-2<\Re z\leq -1, z \ne -1 \}$.
If $q$ is sufficiently close to 1 then $[z+1]_{q} \ne 0$ on $K$. Therefore as 
$q\to1-0$, 
 \begin{equation*}
   \Gamma(z+1;q)=\frac{\Gamma(z+2;q)}{[z+1]_{q}}
 \end{equation*}  
uniformly converges to
  \begin{equation*}
    \frac{\Gamma(z+2)}{z+1} = \Gamma(z+1)
  \end{equation*}  
on $K$. \par 
We assume that 
  \begin{equation*}
    G_{n-1}(z+1;q) \to G_{n-1}(z+1) \quad \mbox{as} \quad q\to 1-0
  \end{equation*}
and that the convergence is uniform on $K$\par
From (2.7), we see that, if $q$ is suffciently close to 1, $G_{n-1}(z+1;q)$
has no poles and no zeros on $K$, neither has $G_{n-1}(z+1)$ from (2.5).
Therefore, as $q \to 1-0$, 
  \begin{equation*}
    G_{n}(z+1;q)=\frac{G_{n}(z+2;q)}{G_{n-1}(z+1;q)} 
  \end{equation*}
uniformly converges to
  \begin{equation*}
    \frac{G_{n}(z+2)}{G_{n-1}(z+1)} = G_{n}(z+1) 
  \end{equation*}
on $K$.\par
Repeating this procedure, we can verify, for any $n$, $G_{n}(z+1;q)$ converge
to $G_{n}(z+1)$ in a compact set in the domain 
$\{-3<\Re z \leq -2, z\ne -2\}$, $\{-4<\Re z \leq -3, z\ne -3\}, \cdots $.
Thus the claim of the theorem is proved.
\end{pf}


\subsection{Asymptotic expansion of $G_{n}(z+1)$}
Let us call the expression (4.8) the Euler-MacLaurin expansion of 
$G_{n}(z+1).$ We should note that (4.8) is valid for 
$z \in {\bold C}\backslash\{\Re z \leq -1\}.$  
We show that it gives an asymptotic expansion of $G_{n}(z+1)$
as $|z|\to \infty,$ {\it i.e. the higher Stirling formula}.\par 

\begin{thm}
Let $0<\delta<\pi$, then
  \begin{align*} 
    & \log G_{n}(z+1)\\[4pt]
    &  \qquad \sim  \left\{ \binom{z+1}{n}+ \sum_{r=1}^{n}
      \frac{B_{r}}{r!} \left( - \frac{d}{dz} \right)^{r-1}
      \binom{z}{n-1} \right\} \log(z+1) \\[8pt]
    & \qquad - \sum_{r=1}^{n}\left\{
      \left( - \frac{d}{dz} \right)^{r-1}
      \binom{z}{n-1}\right\}\times
      \frac{1}{r!r} \left\{(z+1)^{r}-1\right\}\\[8pt] 
    & \qquad -  \sum_{j=0}^{n-1}G_{n,j}(z)\left\{
      \zeta'(-j)+\frac{1}{(j+1)^{2}}\right\} 
      + \sum_{r=1}^{\infty} \frac{B_{2r}}{(2r)!} F_{n,2r-1}(z)
  \end{align*} 
\vspace{8pt}

as $|z|\to\infty$ in the sector $\{z \in {\bold C}||arg z|<\pi - 
\delta\}$.
\end{thm}

\begin{pf}
Straightforward calculation shows that
{\allowdisplaybreaks
  \begin{align*}
    F_{n,r-1}(z)&=\sum_{l=1}^{r-1}\binom{r-1}{l}
      \left\{\left(\frac{d}{dt}\right)^{r-1-l}
      \binom{-t}{n-1}\right\}\Bigg\vert_{t=1}
      \left\{\left(\frac{d}{dt}\right)^{l}
      \log(z+t)\right\}\Bigg\vert_{t=1}\\[8pt]
    & = \sum_{l=1}^{r-1}\left(
      \sum_{k=0}^{n-1}\sideset{_{n-1}}{_k}S (-1)^{k}[k]_{r-1-l}\right)
      \frac{(-1)^{l-1}(l-1)!}{(z+1)^{l}}.
  \end{align*}}
\vspace{8pt}

Thus, if $r>n$, then $F_{n,r-1}(z) = O(z^{-r+n})$ as $|z|\to
\infty$.  So we can see that
  \begin{equation}
    \frac{F_{n,2r-1}(z)}{F_{n,2r-3}(z)} = O(z^{-2}) \quad \mbox{as}
    \quad |z| \to \infty.
  \label{eqn;24}\end{equation}

Furthermore we can see that 

{\allowdisplaybreaks
  \begin{align*}
    & R_{n,2m}(z) = \frac{-1}{(2m)!}\int_{1}^{\infty}
        \overline{B}_{2m}(t)\left(\frac{d}{dt}\right)^{2m}
        \left\{\binom{-t}{n-1} \log \left(
        \frac{z+t}{z+1}\right)\right\}dt\\[8pt]
    & \qquad = \frac{-1}{(2m)!}\sum_{j=1}^{n-1}(-1)^{j}
      \sideset{_{n-1}}{_j}S  \sum_{l=0}^{j}[j]_{l}(m-l-1)!
      \int_{1}^{\infty}{\overline B}_{2m}(t)
      \frac{t^{j}}{(t+z)^{2m-j}}dt.
  \end{align*}}
Noting that
  \begin{equation*}
    \frac{1}{|z+t|} < \frac{1}{|t||\sin \delta|}
  \end{equation*}
in the sector $\{z\in{\bold C}||arg z|<\pi - \delta\}$ 
and that $|B_{2m}(t)| \leq |B_{2m}|$ for $0\leq t \leq 1$,
we have 
{\allowdisplaybreaks
  \begin{align*}
    & \Bigg\vert \int_{1}^{\infty} {\overline B}_{2m}(t)
      \left(\frac{d}{dt}\right)^{2m} \left\{\binom{-t}{n-1}
      \log \left(\frac{z+t}{z+1}\right)\right\}dt\Bigg\vert\\[8pt]
    & \qquad \leq \frac{|B_{2m}|}{(2m)!}\sum_{j=1}^{n-1}\sum_{l=0}^{j}
      \sideset{_{n-1}}{_j}S \binom{n}{l}
      \frac{[j]_{l}(2m-l-1)!}{|\sin \delta|^{2m-j}}
      \int_{1}^{\infty}\frac{dt}{t^{2m-j}}. 
  \end{align*}}

Hence, 
  \begin{equation*}
    |R_{n,2m}(z)|=O(z^{-2m-1+n})=o(F_{n,2m-1}(z)) 
  \end{equation*}
as $|z|\to\infty$ in the sector.
\end{pf}

Let us exhibit some examples of the higher Stirling formula. 
In the case that $n=1$, we obtain
  \begin{eqnarray*}
    & &\log G_{1}(z+1)=\log \Gamma(z+1)\\ 
    & &\qquad\sim
      \left(z+\frac{1}{2}\right)\log(z+1) - (z+1) - \zeta'(0)\\ 
    & &\qquad + \sum_{r=1}^{\infty} \frac{B_{2r}}{[2r]_{2}}
    \frac{1}{(z+1)^{2r-1}}.
  \end{eqnarray*}    
This is the Stirling formula since $\zeta'(0)=-\frac{1}{2}\log
(2\pi)$.\par 
Furthermore, in the case that $n=2$, we obtain 
  \begin{align*}
     & \log G_{2}(z+1)\\[4pt] 
     & \quad \sim
       \left(\frac{z^{2}}{2}-\frac{1}{12}\right) \log(z+1)
       -\frac{3}{4}z^{2}-\frac{z}{2}+\frac{1}{4}\\[8pt] 
     & \quad - z \zeta'(0)+ \zeta'(-1) \\[4pt]
     & \quad -\frac{1}{12}\frac{1}{z+1}
       + \sum_{r=2}^{\infty}\frac{B_{2r}}{[2r]_{3}}
       \frac{1}{(z+1)^{2r-1}}(z-2r+1),
  \end{align*}
\vspace{8pt}

\noindent which coincides with the formula (2.3). In the case that $n=3,4,  $
and $5$, we have the following results.

\begin{prop}
The higher Stiring formula for $n=3,4$ and $5$ are as follows: 
{\allowdisplaybreaks
  \begin{align*}
   & \log G_{3}(z+1)\\[4pt]
   & \quad \sim
     \left(\frac{z^{3}}{6}-\frac{z^{2}}{4}+\frac{1}{24}\right)\log(z+1)
     - \frac{11}{36}z^{3}+\frac{5}{24}z^{2}
     + \frac{z}{3} - \frac{13}{72}\\[8pt]
   &  \quad - \frac{z^{2}-z}{2} \zeta'(0)  
     + \frac{2z-1}{2} \zeta'(-1) -\frac{1}{2}\zeta'(-2)\\[8pt]
   & \quad  +\frac{1}{12}\frac{1}{z+1}
     +  \sum_{r=2}^{\infty}
      \frac{B_{2r}}{[2r]_{4}}\frac{1}{(z+1)^{2r-1}}
      \left\{z^{2} - (6r-11)z+(4r^2-16r+16)\right\}.\\
&\\
  & \log G_{4}(z+1)\\[4pt]
  & \quad \sim
    \left(\frac{z^4}{24}-\frac{z^3}{6}+\frac{z^2}{6}
    -\frac{19}{720}\right)\log(z+1)\\[8pt]
  & \quad - \frac{4}{72}z^4+\frac{2}{9}z^3
    + \frac{z^2}{8}-\frac{11}{36}z+\frac{31}{144}\\[8pt]
  & \quad - \frac{z^3 - 3z^2 +2z}{6}\zeta'(0) 
    + \frac{3z^2-6z+2}{6}\zeta'(-1) 
    -\frac{z-1}{2}\zeta'(-2)+\frac{1}{6}\zeta'(-3)\\[8pt]
  & \quad - \frac{1}{12}\frac{1}{z+1}+
    \frac{1}{720}\frac{1}{(z+1)^3}\left(
    6z^2+\frac{13}{2}z+\frac{5}{2}\right)\\[8pt]
  & \quad + \sum_{r=3}^{\infty}\frac{B_{2r}}{[2r]_{5}}
    \frac{1}{(z+1)^{2r-1}}\left\{
    z^3-(12r-27)z^2+(20r^2-94r+111)z\right.\\[4pt]
  & \left.\qquad \qquad
    -(8r^3-56r^2+134r-109)
    \right\}.\\
&\\
  & \log G_{5}(z+1)\\[4pt]
  & \quad \sim \left(
    \frac{z^{5}}{120}-\frac{z^{4}}{16}+\frac{11}{72}z^{3}
    -\frac{z^{3}}{8}+\frac{3}{160}\right)\log(z+1)\\[8pt]
  & \quad -\frac{137}{7200}z^{5}+\frac{39}{320}z^{4}-\frac{461}{2160}z^{3}
    +\frac{z^{2}}{1440}-\frac{323}{1440}z+\frac{5639}{43200}\\[8pt]
  &  \quad  -\frac{z^{4}-6z^{3}+11z^{2}-6z}{24}\zeta'(0)
    +\frac{4z^{3}-18z^{2}+22z-6}{24}\zeta'(-1)\\[8pt]
  &  \quad - \frac{6z^{2}-18z+11}{24}\zeta'(-2)
    +\frac{2z-3}{12}\zeta'(-3)-\frac{1}{24}\zeta'(-4)\\[8pt]
  &  \quad + \frac{1}{12}\frac{1}{z+1}
    -\frac{1}{720}\frac{1}{(z+1)^{3}}\left(\frac{35}{4}z^{2}
    +\frac{45}{4}z+\frac{9}{2}\right)\\[8pt]
  &  \quad + \sum_{r=3}^{\infty}\frac{B_{2r}}{[2r]_{6}}
    \frac{1}{(z+1)^{2r-1}}
    \left\{z^{4}-\left(20r-54\right)z^{3}
    +\left(70r^{2}-375r+506\right)z^{2}\right.\\[8pt]
  &  \quad \qquad - \left(\frac{200}{3}r^{3}-540r^{2}
    +\frac{4420}{3}r-1354\right)z\\[4pt]
  &  \quad \qquad \qquad \left. +16r^{4}-\frac{536}{3}r^{3}+754r^{2}
    -\frac{4279}{3}r+1021\right\}.
  \end{align*}}
\end{prop}


\section{The Weierstrass product repesentation for the $G_{n}(z+1)$}
By calculating the formula (\ref{eqn;22}) in the case of $m=n+1$, 
we derive the Weierstrass product representation for the multiple 
gamma functions. Main theorem of this section is the following:

\begin{thm}
For $n \in {\bold N}$, we have 
  \begin{equation*}
    G_{n}(z+1) 
      = \exp \left(F_{n}(z) \right)
         \prod_{k=1}^{\infty}\left\{\left(
           1+\frac{z}{k}\right)^{-\binom{-k}{n-1}}
           \exp \left( \Phi_{n}(z,k) \right)
         \right\},
  \end{equation*}         
where
{\allowdisplaybreaks
\begin{align*}
  &  F_{n}(z):= \sum_{j=0}^{n-1}G_{n,j}(z)Q_{j}(z)
    + \sum_{r=0}^{n-2}\left[
      \frac{1}{r!}\left(\frac{\partial}{\partial u}\right)^{r}
      \binom{z-u}{n-1}\right]_{u=0}^{u=z} \times \zeta'(-r)\\
  &  \quad \qquad - \int_{0}^{z} \binom{z-u}{n-1} du \times \gamma,\\[8pt]  
  & \Phi_{n}(z,k) := \frac{1}{(n-1)!} \sum_{\mu=-1}^{n-2}
    \left\{\sum_{r = \mu + 1}^{n-1}
    \frac{\sideset{_{n-1}}{_r}S}{r-\mu}z^{r-\mu}\right\}
    (-1)^{\mu+1}k^{\mu},\\[8pt]
  & Q_{j}(z) := P_{j}(z+1) -\sum_{r=0}^{j}\binom{j}{r}z^{r}P_{j-r}(1)\\[4pt]
  & \qquad \quad + \frac{1}{j+1}\sum_{r=1}^{j+1}\binom{j+1}{r}
    B_{j+1-r}(z)\sum_{l=1}^{r}\frac{(-1)^{l-1}z^{l}}{l},\\[8pt]
  & P_{j}(x) := \sum_{r=0}^{j+1}\frac{B_{r}}{r!}
    \varphi_{j.r} x^{j-r+1},\\[8pt]
  & \varphi_{j,r} := \left(\frac{d}{dt}\right)^{r}\left\{
    \frac{t^{j+1}}{j+1}\log t - \frac{t^{j+1}}{(j+1)^{2}}
    \right\}\Bigg\vert_{t=1}.
  \end{align*}}
\label{thm;31}\end{thm} 

Proof of this theorem will be carried out through the sections 5.1 $\sim$ 
5.4, and some examples of this representation will be discussed 
in the section 5.5.\par


\subsection{Rewriting the Euler-MacLaurin expansion of $G_{n}(z+1)$}
In this section we prove the following proposition.

\begin{prop}
  \begin{equation}
    \log G_{n}(z+1)=\sum_{j=0}^{n-1}G_{n,j}(z)K_{j}(z)
  \label{eqn;31}\end{equation}
where
  \begin{align*}
     K_{j}(z) & := \frac{B_{j+1}(z+1)}{j+1}\log(z+1) 
       - \zeta'(-j)+ P_{j}(z+1)\\[4pt]
       & + \sum_{k=1}^{\infty} \left[
         P_{j}(z+k+1)-P_{j}(z+k) + 
         \frac{B_{j+1}(z+k+1)}{j+1}
         \log \left(
         \frac{z+k+1}{z+k}\right)\right].
  \end{align*}
Furthermore the infinite sum of the last term is absolutely
convergent.
\label{prop;31}\end{prop}

\begin{pf}
From the definiftion of $G_{n,r}(z)$, we have
  \begin{equation}
    \left(-\frac{d}{dt}\right)^{r-1}\binom{z}{n-1}
    = \left(\frac{d}{du}\right)^{r-1}\binom{z-u}{n-1}\Bigg\vert_{u=0}
    = (r-1)!G_{n,r-1}(z),
  \label{eqn;32}\end{equation}  
and hence
  \begin{align}
     & \sum_{r=1}^{n} \left\{ \left(-\frac{d}{dz}\right)^{r-1}
       \binom{z}{n-1}\right\}\Bigg\vert_{z=0}\times \frac{1}{r!r}
       \left\{(z+1)^{r}-1\right\}\label{eqn;33}\\[4pt]
     & \qquad \qquad =\sum_{j=0}^{n-1}G_{n,j}(z)
       \frac{(z+1)^{j+1}-1}{(j+1)^{2}}.\nonumber
  \end{align}
Since 
  \begin{equation*}
    \binom{-t}{n-1} = \binom{z-(z+t)}{n-1}
    = \sum_{j=0}^{n-1}G_{n,j}(z)(z+t)^{j},
  \end{equation*}
we have, for $1 \leq r \leq n-1$
  \begin{equation}
    F_{n,r-1}(z)
      = \sum_{j=0}^{n-1}G_{n,j}(z)(z+1)^{j+1-r}
      \varphi_{j,r}.
  \label{eqn;34}\end{equation}
In order to calculate the coefficient of $\log(z+1)$, we 
use the following lemma.

\begin{lem}
  \begin{equation*}
    \binom{z+1}{n} + \sum_{j=0}^{n-1} G_{n,j}(z)
      \frac{B_{j+1}}{j+1}
    = \sum_{j=0}^{n-1}G_{n,j}(z)\frac{B_{j+1}(z+1)}{j+1}.
  \end{equation*}  
\label{lem;31}\end{lem}

\begin{pf}
Since the both sides are polynomials in $z$, it is sufficient to 
prove that the formula holds for a sufficiently large, arbitrary 
integer $N$. Noting that
  \begin{equation*}
    \frac{B_{j+1}(N+1)-B_{j+1}}{j+1} = 
      \begin{cases}
        \sum_{l=1}^{N} l^j \quad (j>0)\\
        N+1 \quad (j=0),
      \end{cases}
  \end{equation*}    
we have
{\allowdisplaybreaks
  \begin{align*}
    & \sum_{j=0}^{n-1}G_{n,j}(N)\frac{B_{j+1}(N+1)}{j+1}
      - \sum_{j=0}^{n-1}G_{n,j}(N)\frac{B_{j+1}}{j+1}\\[4pt] 
    & \quad = \sum_{l=1}^{N} \sum_{j=0}^{n-1}
      G_{n,j}(N)l^{j} + G_{n,0}(N)\\[4pt]
    & \quad = \sum_{l=1}^{N}\binom{N-l}{n-1} + \binom{N}{n-1}\\[4pt]
    & \quad = \binom{N}{n}+\binom{N}{n-1}\\[4pt]
    & \quad = \binom{N+1}{n}.
 \end{align*}}  
\end{pf}

From (\ref{eqn;32}) and Lemma \ref{lem;31},
it follows that 
{\allowdisplaybreaks
  \begin{align}
    & \binom{z+1}{n-1} + \sum_{r=1}^{n} \frac{B_{r}}{r!}
      \left(-\frac{d}{dt}\right)^{r-1}\binom{z}{n-1}\label{eqn;35}\\[8pt]
    & \quad = \binom{z+1}{n-1} + \sum_{j=0}^{n-1}
      G_{n,j}(z) \frac{B_{j+1}}{j+1}\nonumber\\[8pt]
    & \quad = \sum_{j=0}^{n-1}G_{n,j}(z)
      \frac{B_{j+1}(z+1)}{j+1}.\nonumber
  \end{align}}    

Next we calculate $R_{n,n+1}(z)$. From the definition of
$G_{n,r}(z)$, we have
{\allowdisplaybreaks
  \begin{align}
    &  R_{n,n+1}(z)\label{eqn;36}\\[4pt] 
    & = \frac{(-1)^{n}}{(n+1)!}
      \int_{1}^{\infty}\overline{B}_{n+1}(t)
      \left(\frac{d}{dt}\right)^{n+1}
      \left\{\binom{z-(z+t)}{n-1}
      \log \left( \frac{z+t}{z+1} \right)
      \right\}dt \nonumber\\[8pt]
    & + \sum_{j=0}^{n-1}G_{n,j}(z)\left[
      \frac{(-1)^{n}}{(n+1)!} \int_{1}^{\infty}
      \overline{B}_{n+1}(t)\left(
      \frac{d}{dt}\right)^{n+1} \left\{
      (z+t)^{j}\log \left(
      \frac{z+t}{z+1}\right)\right\}\right]dt\nonumber\\[8pt]
    & = \sum_{j=0}^{n-1}G_{n,j}(z)\left[
      \sum_{k=1}^{\infty}\frac{(-1)^{n}}{(n+1)!}
      \int_{0}^{1} B_{n+1}(t)\left(\frac{d}{dt}\right)^{n+1}
      \left\{(z+t+k)^{j}\log\left(
      \frac{z+t+k}{z+1}\right)\right\}dt\right].\nonumber
 \end{align}}    

Here we have, for $j \leq n-1$, 
  \begin{align}
     & \int_{0}^{1} B_{n+1}(t) \left(\frac{d}{dt}\right)^{n+1}
      \left\{ (z+t+k)^{j} \log \left( \frac{z+t+k}{z+1} \right)
      \right\}dt \label{eqn;37}\\[4pt]
      & \quad \qquad = O(k^{j-n-1}) \quad \mbox{as} \quad
      k \to \infty, \nonumber
  \end{align}    
so that the infinite sum in (\ref{eqn;36}) converges absolutely because 
$j-n-1 \leq -2$.\par
By means of the Euler-MacLaurin summation formula, we have
{\allowdisplaybreaks
  \begin{align}
    & \frac{(-1)^{n}}{(n+1)!}\int_{0}^{1}
      B_{n+1}(t)\left(\frac{d}{dt}\right)^{n+1}
      \left\{(z+t+k)^{j}\log\left(
      \frac{z+t+k}{z+1}\right)\right\}dt \label{eqn;38}\\[8pt]
    & = \left\{(z+k)^{j}+\frac{(z+k)^{j+1}}{j+1} +
      \sum_{r=1}^{j+1}\frac{B_{r}}{r!}[j]_{r-1}(z+k)^{j+1-r}
      \right\}\log(z+k)\nonumber\\[8pt]
    & - \left\{\frac{(z+k+1)^{j+1}}{j+1}
      + \sum_{r=1}^{j+1}\frac{B_{r}}{r!}[j]_{r-1}
      (z+k+1)^{j+1-r}\right\}\log(z+k+1)\nonumber\\[8pt]
    & + \left[ -(z+k)^{j} + \frac{(z+k+1)^{j+1}}{j+1}
      - \frac{(z+k)^{j+1}}{j+1}\right.\nonumber\\[8pt]
    &  \quad \left.- \sum_{r=1}^{j+1}\frac{B_{r}}{r!}[j]_{r-1}
      \left\{(z+k)^{j+1-r}-(z+k+1)^{j+1-r}\right\}
      \right]\log(z+1)\nonumber\\[8pt]
    & - \sum_{r=0}^{n+1}\frac{B_{r}}{r!}\varphi_{j,r}
      \left\{(z+k+1)^{j+1-r}-(z+k)^{j+1-r}\right\}.
      \nonumber
  \end{align}}

The coefficient of $\log(z+1)$ of the above formula vanishes because of 
the following lemma.

\begin{lem}
  \begin{align}
     &\frac{B_{j+1}(z+k+1)}{j+1} 
      = \frac{(z+k+1)^{j+1}}{j+1} 
      + \sum_{r=1}^{j+1} \
      \frac{B_{r}}{r!}[j]_{r-1}(z+k+1)^{j+1-r}\nonumber\\[4pt]
     & \quad = (z+k)^{j} + \frac{(z+k)^{j+1}}{j+1}
      + \sum_{r=1}^{j+1}
      \frac{B_{r}}{r!}[j]_{r-1}(z+k)^{j+1-r}. \nonumber
  \end{align}  
\label{lem;32}\end{lem}

\begin{pf}
The first equality follows from the identity
  \begin{equation*}
    B_{j+1}(z+k+1) = \sum_{r=0}^{j+1}
      \binom{j+1}{r}(z+k+1)^{j+1-r}B_{r}.
  \end{equation*} 
The second equality holds because of the Euler-MacLaurin summation
formula for $(z+t+k)^{j}$ on $[0,1]$.
\end{pf}   

Applying Lemma \ref{lem;32} to (\ref{eqn;38}),
we obtain 
  \begin{align*}
    & R_{n,n+1}(z) = - \sum_{j=0}^{n-1}G_{n.j}(z)
      \sum_{k=1}^{\infty} \left[
      \frac{B_{j+1}(z+k+1)}{j+1}
      \log\left(
      \frac{z+k+1}{z+k}\right)\right.\\[8pt]
    & \qquad \left.+ \sum_{r=0}^{n+1}\frac{B_{r}}{r!}
      \varphi_{j,r} \left\{
      (z+k+1)^{j+1-r}-(z+k)^{j+1-r}\right\}\right].
  \end{align*}    
If $r\geq j+2$, then $(z+k+1)^{j+1-r}-(z+k)^{j+1-r}$ decreases
more rapidly than $k^{-2}$ as $k\to\infty$. Thus we have
{\allowdisplaybreaks 
  \begin{align}
   & R_{n,n+1}(z) = -\sum_{j=0}^{n-1} G_{n,j}(z)
     \left[ - \sum_{r=j+2}^{n+1} \frac{B_{r}}{r!}
     (z+1)^{j+1-r} \varphi_{j,r}\right.\label{eqn;39}\\[8pt]
   & \left. + \sum_{k=1}^{\infty} \left\{
     \frac{B_{j+1}(z+k+1)}{j+1}\log\left(
     \frac{z+k+1}{z+k}\right)
     + P_{j}(z+k+1)-P_{j}(z+k)\right\}
     \right]\nonumber,
   \end{align}}
and the infinite sum is absolutely convergent.\par
Substituting (\ref{eqn;33}), (\ref{eqn;34}), (\ref{eqn;35}),
(\ref{eqn;39}) to (\ref{eqn;22}), and noting that
  \begin{equation*}
    -\frac{1}{(j+1)^{2}}
      +\sum_{r=1}^{n+1}\frac{B_{r}}{r!}(z+1)^{j+1-r}
      -\sum_{r=j+2}^{n+1}\frac{B_{r}}{r!}(z+1)^{j+1-r}
      = P_{j}(z+1),
  \end{equation*} 
we obtain (5.1).
\end{pf}

The next lemma which will be useful in the section 5.3 is 
deduced from the above proof. 

\begin{lem}
  \begin{align*}
    &  \frac{B_{j+1}(z+k+1)}{j+1}\log\left(
      \frac{z+k+1}{z+k}\right)+P_{j}(z+k+1)-P_{j}(z+k)\\[4pt]
    & \qquad = O(k^{-2}) \quad \mbox{as} \quad k \to \infty.  
  \end{align*}
In other words, the polynomial $P_{j}(z+k+1)-P_{j}(z+k)$ is a convergent
factor such that the infinite sum in (\ref{eqn;39}) converges absolutely.
\label{lem;33}\end{lem}


\subsection{An infinite product representation for the $\zeta'(-j)$}
We derive an infinite product representation for $\zeta'(-j)$ in the same way
as in the section 3.1. 

\begin{prop}
  \begin{equation*}
    \exp(\zeta'(-j)) = \exp(P_{j}(1))
     \prod_{k=1}^{\infty}\left\{
      \left(1+\frac{1}{k}\right)^{\frac{B_{j+1}(k+1)}{j+1}}
      \exp\left(P_{j}(k+1)-P_{j}(k)\right)\right\}
  \end{equation*}
and the infinite product converges absolutely.
\label{thm;32}\end{prop}

\begin{pf}
From the proof of Lemma \ref{lem;23}, we have
{\allowdisplaybreaks
  \begin{align*}
     & \zeta'(-j) = - \frac{1}{(j+1)^{2}}+
       \sum_{r=1}^{j+2}\frac{B_{r}}{r!}\varphi_{j,r}\\[4pt]
     & \qquad  \qquad \qquad - \frac{(-1)^{j+1}}{(j+2)!}\int_{1}^{\infty}
       \overline{B}_{j+2}(t)\left(\frac{d}{dt}\right)^{j+2}
       \{ t^{j} \log t \}dt\\[4pt]
     & \quad = P_{j}(1)+\frac{B_{j+2}}{(j+2)!}\varphi_{j,j+2}\\[4pt]
     & \qquad \qquad \qquad - \frac{(-1)^{j+1}}{(j+2)!}\int_{1}^{\infty}
       \overline{B}_{j+2}(t)\left(
       \frac{d}{dt}\right)^{j+2}\{t^{j}\log t\}dt.
  \end{align*}}
Furthermore we can see by the same way as (\ref{eqn;37}) that
  \begin{equation}
    \int_{1}^{\infty}\overline{B}_{j+2}(t)\left(
    \frac{d}{dt}\right)^{j+2}\{t^{j}\log t\} dt
    = O(k^{-2}) \quad \mbox{as}\quad k\to\infty,
  \label{eqn;310}\end{equation}
Hence we have
  \begin{align*}
    & \frac{(-1)^{j+1}}{(j+2)!}\int_{1}^{\infty}
    \overline{B}_{j+2}(t)\left(\frac{d}{dt}\right)^{j+2}
    \{t^{j} \log t\}dt\\[8pt]
    & \quad = \frac{(-1)^{j+1}}{(j+2)!} 
      \sum_{k=1}^{\infty} \int_{0}^{1}B_{j+2}(t)
      \left(\frac{d}{dt}\right)^{j+2} \{(t+k)^{j} \log (t+k) \}dt,
 \end{align*}
and the infinite sum of this formula is absolutely convergent.
In the same fashion as in the previous section, we have 
  \begin{align}
    & \frac{(-1)^{j+1}}{(j+2)!}\int_{0}^{1}
      B_{j+2}(t)
      \left(\frac{d}{dt}\right)^{j+2}
      \{(t+k)^{j} \log (t+k) \}dt \label{eqn;311}\\[8pt]
    & \qquad = - \frac{B_{j+1}(k+1)}{j+1} \log
      \left(1+\frac{1}{k}\right)+P_{j}(k)-P_{j}(k+1)\nonumber\\[8pt]
    & \qquad + \frac{B_{j+2}}{(j+2)!}
      \left(\frac{1}{k}-\frac{1}{k+1}\right)
      \varphi_{j,j+2}.\nonumber
  \end{align}
From (\ref{eqn;310}) and (\ref{eqn;311}), it follows that
  \begin{align}
    & \frac{B_{j+1}(k+1)}{j+1}\log\left(1+\frac{1}{k}\right)
     + P_{j}(k+1)-P_{j}(k) \label{eqn;312}\\[4pt]
    & \qquad =O(k^{-2}) \quad \mbox{as} \quad 
      k\to\infty. \nonumber
  \end{align}
Hence we have 
  \begin{equation*}
    \zeta'(-j) = P_{j}(1) + \sum_{k=1}^{\infty}
      \log \left[ \left(1 + \frac{1}{k} \right)^{\frac{B_{j+1}(k+1)}{j+1}}
      \exp \left\{ P_{j}(k+1)-P_{j}(k) \right\} \right],
  \end{equation*}
and the infinite sum is absolutely convergent.
\end{pf}

\begin{rem}
In a similar fashion to Proposition 5.6, an infinite product
representaion of $\zeta'(-j,z) \quad (:=\frac{d}{ds}\zeta(s,z)|_{s=-j})$ 
can be given as follows:
  \begin{align*}
    & \exp(\zeta'(-j,z))
      =\exp\left(-\frac{z^{j+1}}{j+1}\log z+P_{j}(z)\right) \\[4pt]
    & \qquad \times  \prod_{k=0}^{\infty}\left[\left(
        \frac{z+k+1}{z+k}\right)^{\frac{B_{j+1}(z+k+1)}{j+1}}
        \exp \left(P_{j}(z+k+1)-P_{j}(z+k)\right)
      \right].  
  \end{align*}
\end{rem}


\subsection{Good representation for $K_{j}(z)$}
In this section, we give a ``good'' representation for $K_{j}(z)$.
The following two lemmas are useful in subsequent arguments.

\begin{lem}
Let $k$ be a positive integer and $k\to \infty$, then we have
  \begin{align*}
    &(1) \quad k^{r}\log
      \left(1+\frac{z}{k}\right)
      +\sum_{l=1}^{r+1}\frac{(-1)^{l}z^{l}}{l}k^{r-l}
      = O(k^{-2})\\[8pt]
    & (2) \quad \frac{B_{j+1}(z+k+1)}{j+1}
      \log\left(1+\frac{z}{k}\right)
      + \frac{1}{j+1} \sum_{r=0}^{j+1}B_{j+1-r}(1)
      \sum_{l=1}^{r+1}\frac{(-1)^{l}}{l}(z+k)^{l-r}\\[4pt]
    & \qquad = O(k^{-2}).
  \end{align*}    
\label{lem;34}\end{lem}

\begin{pf}
(1) is clear. (2) follows from the identity
  \begin{equation*}
    B_{j+1}(z+k+1)=\sum_{r=0}^{j+1}
      B_{j+1-r}(1)(z+k)^{r}.
  \end{equation*}
\end{pf}

\begin{lem}
There exists a unique polynomial $A(k,z)$ such that
  \begin{equation*}
    \frac{B_{j+1}(z+k+1)}{j+1} \log\left(
      1+\frac{z}{k}\right) + A(k,z)
      = O(k^{-2}).
  \end{equation*}
$A(k,z)$ is equal to $P_{j}(z+k+1)-P_{j}(z+k)$.
\label{lem;35}\end{lem}
    
\begin{pf}
By integrating the both side of 
  \begin{equation*}
    B_{j+1}(z+1) = \sum_{r=0}^{j+1}B_{j+1-r}(1)z^{r}
  \end{equation*}
from $-1$ to $0$, we get the formula
  \begin{equation}
    \frac{1}{j+1}\sum_{r=0}^{j+1}
    \binom{j+1}{r}B_{j+1-r}(1)
    \frac{(-1)^{r}}{r+1} = 0.
  \label{eqn;313}\end{equation}  
From Lemma \ref{lem;34} (2) and (\ref{eqn;313}), a polynomial
$A(k,z)$ satisfying
  \begin{equation}
    \frac{B_{j+1}(z+k+1)}{j+1}\log \left(1+\frac{z}{k}\right)
    + A(k,z) = O(k^{-2})
    \quad \mbox{as} \quad k\to \infty,
  \label{eqn;314}\end{equation}
is uniquely determined and  the polynomial $P_{j}(z+k+1)-P_{j}(z+k)$ 
satisfies (\ref{eqn;314}) by Lemma \ref{lem;33}. Therefore,
$A(k,z)=P_{j}(z+k+1)-P_{j}(z+k)$.
\end{pf}

Using these lemmas, we prove the following proposition.

\begin{prop}
Let $k$ be a positive integer and define
{\allowdisplaybreaks 
  \begin{align*}
    & (I)_{k}:= \frac{1}{j+1}\left\{
      \sum_{r=0}^{j+1}\binom{j+1}{r}B_{j+1-r}(z)
      \sum_{l=1}^{r+1}\frac{(-1)^{l-1}z^{l}}{l}k^{l-r}
      \right\},\\[8pt]
    & (II)_{k}:= - (I)_{k+1},\\[8pt]
    & (III)_{k}:= \sum_{r=0}^{j}z^{j-r}
        \left\{P_{r}(k+1)-P_{r}(k)\right\}
        - \frac{1}{k}\frac{z^{j+1}}{j+1},\\[8pt]
and\\[8pt]
    &(IV)_{k} := \sum_{r=0}^{j}\binom{j}{r}z^{j-r}
        \sum_{l=1}^{r+1} \frac{(-1)^{l-1}z^{l}}{l}k^{r-l}.
  \end{align*}}      
Then we have 
{\allowdisplaybreaks
  \begin{align}
    & \frac{B_{j+1}(z+k+1)}{j+1}
      \log \left(\frac{z+k+1}{z+k}\right)
      + P_{j}(z+k+1)-P_{j}(z+k)\label{eqn;315}\\[8pt]
    & \quad =  \left\{\frac{B_{j+1}(z+k)}{j+1}
      \log \left(\frac{k}{z+k}\right)+(I)_{k}
      \right\}\nonumber\\[8pt]
    & \quad + \left\{\frac{B_{j+1}(z+k+1)}{j+1}
      \log \left(\frac{z+k+1}{k+1}\right)+(II)_{k}
      \right\}\nonumber\\[8pt]
    & \quad + \left\{\frac{B_{j+1}(z+k+1)}{j+1}
      \log \left(\frac{k+1}{k}\right)+(III)_{k}
      \right\}\nonumber\\[8pt]
    & \quad+ \left\{-(z+k)^{j}
      \log\left(1+\frac{z}{k}\right)+(IV)_{k}
      \right\}\nonumber\\[8pt]
    & \quad- \frac{1}{j+1}\left\{
      \sum_{r=0}^{j+1}\frac{(-1)^{r}z^{r+1}}{r+1}
      B_{j+1-r}(z)\right\}
      \left(\frac{1}{k}-\frac{1}{k+1}\right)\nonumber
  \end{align}}    
\vspace{8pt}

Furthermore each term in the right hand side decreases like $O(k^{-2})$ as
$k\to \infty$.
\label{prop;32}\end{prop}

\begin{pf} 
Making use of the identity,
  \begin{equation*}
    \frac{B_{j+1}(z+k)-B_{j+1}(z+k+1)}{j+1}
      = -(z+k)^{j},
  \end{equation*}
we have
{\allowdisplaybreaks
  \begin{align}
    & \frac{B_{j+1}(z+k+1)}{j+1}
       \log\left(\frac{z+k+1}{z+k}\right)\label{eqn;316}\\[8pt]
    & = \frac{B_{j+1}(z+k)}{j+1}
      \log\left(\frac{k}{z+k}\right)
      + \frac{B_{j+1}(z+k+1)}{j+1}
      \log\left(\frac{z+k+1}{k+1}\right)\nonumber\\[8pt]
    & + \frac{B_{j+1}(z+k+1)}{j+1}
      \log \left( \frac{k+1}{k} \right)
      - (z+k)^{j} \log \left(\frac{z+k}{k}\right).\nonumber
  \end{align}}
\vspace{8pt}

By Lemma \ref{lem;34} (1), (\ref{eqn;312}) and the identity
  \begin{equation*}
    \frac{B_{j+1}(z+k+1)}{j+1}
      = \frac{1}{j+1}\sum_{r=0}^{j+1}
        \binom{j+1}{r}B_{j+1-r}(z)k^{r},
  \end{equation*}
we can see that each term in the right hand side of (\ref{eqn;315})
is $O(k^{-2})$.\par
On the other hand, it can be seen that $(I)_{k}\sim(IV)_{k}$ are 
polynomials of $z$ and that
{\allowdisplaybreaks 
  \begin{align*}
    & (I)_{k}= \frac{1}{j+1} \left\{
        \sum_{r=0}^{j+1}\binom{j+1}{r}
        B_{j+1-r}(z) \frac{(-1)^{r}z^{r+1}}{r+1}
        \right\} \frac{1}{k}\\
    & \qquad  + (\mbox{a polynomial of $k$}),\\[8pt]
    & (II)_{k}= - \frac{1}{j+1} \left\{
      \sum_{r=0}^{j+1}\binom{j+1}{r}
      B_{j+1-r}(z)\frac{(-1)^{r}z^{r+1}}{r+1}
      \right\}\frac{1}{k+1}\\
    & \qquad  + (\mbox{a polynomial of $k$}),\\[8pt]
    & (III)_{k}= - \frac{z^{j+1}}{j+1}\frac{1}{k}
       + (\mbox{a polynomial of $k$}),\\[8pt]
    & (IV)_{k} = \left\{\sum_{r=0}^{j}
        \frac{(-1)^{r}}{j+1} \binom{j}{r}\right\}
        \frac{z^{j+1}}{k}
         + (\mbox{polynomial of $k$})\\[8pt]
    & \qquad = \frac{z^{j+1}}{j+1}\frac{1}{k}
          + (\mbox{a polynomial of $k$}).
  \end{align*}}
Hence, if we put
{\allowdisplaybreaks
  \begin{align*}
    & B(k,z) := (I)_{k}+(II)_{k}+(III)_{k}+(IV)_{k}\\[4pt]
    & \qquad - \frac{1}{j+1} \left\{
      \sum_{r=0}^{j+1}\binom{j+1}{r}
      \frac{(-1)^{r}z^{r+1}}{r+1}B_{j+1-r}(z)\right\}
      \left(\frac{1}{k}-\frac{1}{k+1}\right),
  \end{align*}}
then, $B(k,z)$ is a polynomial of $k$, $z$ and it satisfies
  \begin{equation*}
    \frac{B_{j+1}(z+k+1)}{j+1}
    \log\left(1+\frac{z}{k}\right)
    + B(k,z) = O(k^{-2})
    \quad \mbox{as} \quad k\to\infty.
  \end{equation*}  
By Lemma \ref{eqn;35}, we have
  \begin{equation}
    B(k.z)=P_{j}(z+k+1)-P_{j}(z+k).
  \label{eqn;317}\end{equation}  
By (\ref{eqn;316}) and (\ref{eqn;317}), we can deduce (\ref{eqn;315}).
\end{pf}

\begin{prop}
{\allowdisplaybreaks
  \begin{align*}
    & K_{j}(z) = Q_{j}(z)+ \sum_{r=1}^{j}\binom{j}{r}z^{j-r}
      \zeta'(-r)-\frac{z^{j+1}}{j+1}\gamma\\[4pt]
    & \qquad + \sum_{k=1}^{\infty}\left\{
      -(z+k)^{k}\log \left(1+\frac{z}{k}\right)
      + \sum_{r=0}^{j}\binom{j}{r}z^{j-r}
      \sum_{l=1}^{r+1}\frac{(-1)^{l-1}z^{l}}{l} k^{r-l}
      \right\}
  \end{align*}}
and the infinite sum of this formula converges absolutely.
\label{prop;33}\end{prop}

\begin{pf}
The summation of each term in (\ref{eqn;315}) from $k=1$ to $k=\infty$,
is absolutely convergent because of Proposition 5.9, and
the following caluculation is possible. Namely, 
{\allowdisplaybreaks
  \begin{align}
   & \sum_{k=1}^{\infty}\left\{\frac{B_{j+1}(z+k)}{j+1}
     \log\left(\frac{k}{z+k}\right) + (I)_{k}\right\}
     \label{eqn;318}\\[8pt]
   & + \sum_{k=1}^{\infty}\left\{
     \frac{B_{j+1}(z+k+1)}{j+1}
     \log\left(\frac{z+k+1}{k+1}\right)+(II)_{k}
     \right\}\nonumber\\[8pt]
   & = - \frac{B_{j+1}(z+1)}{j+1}\log (z+1) 
     + \frac{1}{j+1}\sum_{r=0}^{j+1} \binom{j+1}{r}
     B_{j+1-r}(z)\sum_{l=1}^{r+1}\frac{(-1)^{l-1}z^{l}}{l}.
     \nonumber
  \end{align}}
Noting Theorem \ref{thm;32} and 
  \begin{equation*}
    \gamma = \sum_{k=1}^{\infty}\left\{
      \log \left(1+\frac{1}{k}\right)
      - \frac{1}{k}\right\},
  \end{equation*}    
we have
{\allowdisplaybreaks
  \begin{align}
    & \sum_{k=1}^{\infty}\left\{
      \frac{B_{j+1}(z+k+1)}{j+1} 
      \log \left(1+\frac{1}{k}\right) + (III)_{k}
      \right\} \label{eqn;319}\\[4pt]
    & \qquad \quad = \sum_{r=0}^{j}\binom{j}{r}
      z^{j-r}\left\{\zeta'(-r)-P_{r}(1)\right\}
      - \frac{z^{j+1}}{j+1}\gamma.\nonumber
  \end{align}}
By (\ref{eqn;318}) and (\ref{eqn;319}), we obtain
{\allowdisplaybreaks
  \begin{align}
    & \sum_{k=1}^{\infty}\left\{
      \frac{B_{j+1}(z+k+1)}{j+1} 
      \log  \left( \frac{z+k+1}{z+k} \right)
      + P_{j}(z+k+1) - P_{j}(z+k)
      \right\}\label{eqn;320}\\[8pt]
    & \quad = - \frac{B_{j+1}(z+1)}{j+1} \log(z+1)\nonumber\\[8pt]
    & \quad + \frac{1}{j+1}\sum_{r=1}^{j+1}\binom{j+1}{r}
      B_{j+1-r}(z)\sum_{l=1}^{r}\frac{(-1)^{l-1}z^{l}}{l}
      \nonumber\\[8pt]
    & \quad + \sum_{r=0}^{j}\binom{j}{r}z^{j-r}
      \left\{\zeta'(-r) - P_{r}(1)\right\}
      - \frac{z^{j+1}}{j+1}\gamma\nonumber\\[8pt]
    & \quad + \sum_{k=1}^{\infty}\left\{
      - (z+k)^{j} \log \left(1+\frac{z}{k}\right)
      +\sum_{r=0}^{j}\binom{j}{r}z^{j-r}
      \sum_{l=1}^{r+1}\frac{(-1)^{l-1}z^{l}}{l}k^{r-l}
      \right\}.\nonumber
  \end{align}}       
      
The proof is completed by substituting (\ref{eqn;320}) to the
definition of $K_{j}(z)$ in Proposition \ref{prop;31}.
\end{pf}


\subsection{A proof of main theorem}
Using the results in the sections 5.1 $\sim$ 5.4, we prove Theorem 
\ref{thm;31}. By Proposition \ref{prop;33}, we have
{\allowdisplaybreaks
  \begin{align}
    & \log G_{n}(z+1)\label{eqn;321}\\[4pt]
    & \quad = \sum_{j=1}^{n-1}G_{n,j}(z) \left[
        Q_{j}(z)+\sum_{r=0}^{j-1}\binom{j}{r}z^{j-r}
        \zeta'(-r) - \frac{z^{j+1}}{j+1}\gamma\right.
        \nonumber\\[8pt]
    & \left. \quad + \sum_{k=1}^{\infty}\left\{
        -(z+k)^{j} \log \left(1+\frac{z}{k}\right)
        + \sum_{r=0}^{j}\binom{j}{r}z^{j-r}
        \sum_{l=0}^{r+1}\frac{(-1)^{l-1}z^{l}}{l}k^{r-l}
        \right\} \right].\nonumber
   \end{align}}
It is easily seen that
{\allowdisplaybreaks
  \begin{align}
    & \sum_{j=1}^{n-1}G_{n,j}(z)\sum_{r=0}^{j-1}
      \binom{j}{r}z^{j-r}\zeta'(-r),\label{eqn;322}\\[8pt]
    & \qquad = \sum_{r=0}^{n-2} \left\{\sum_{j=0}^{n-1}
      G_{n,j}(z)z^{j-r}\right\}\zeta'(-r)\nonumber,\\[8pt]
    & \qquad = \sum_{r=0}^{n-2}\left[
      \frac{1}{r!}\left(\frac{\partial}{\partial u}\right)
      \binom{z-u}{n-1}\right]_{u=0}^{u=z}\times \zeta'(-r),\nonumber\\[8pt]
\nonumber\\
    & \sum_{j=0}^{n-1} G_{n,j}(z)\frac{z^{j+1}}{j+1}
      = \int_{0}^{z}\binom{z-u}{n-1}du,\label{eqn;323}\\
\nonumber\\
    & \sum_{j=0}^{n-1}G_{n,j}(z)(z+k)^{j}
      = \binom{-k}{n-1}.\label{eqn;324}
  \end{align}}
Substituting (\ref{eqn;322}), (\ref{eqn;323}), (\ref{eqn;324})
to (\ref{eqn;321}), we have
{\allowdisplaybreaks
  \begin{align}
    & \log G_{n}(z+1)\label{eqn;325}\\[4pt]
    & \quad = \sum_{j=0}^{n-1}G_{n,j}(z)Q_{j}(z)
      + \sum_{r=0}^{n-2}\left[
      \frac{1}{r!}\left(\frac{\partial}{\partial u}\right)^{r}
      \binom{z-u}{n-1}\right]_{u=0}^{u=z} \times \zeta'(-r)
      \nonumber\\[8pt]
    & \quad -\int_{0}^{z}\binom{z-u}{n-1}du \times \gamma\nonumber\\[8pt]
    & \quad - \sum_{k=1}^{\infty}\left[
      \binom{-k}{n-1} \log \left(1+\frac{z}{k}\right)
      + \sum_{j=0}^{n-1}G_{n,j}(z)\sum_{r=0}^{j}
      \binom{j}{r} z^{j-r} \sum_{l=1}^{r+1}
      \frac{(-1)^{l-1}z^{l}}{l}k^{r-l}\right].\nonumber
  \end{align}}
Thus, in order to prove Theorem \ref{thm;31}, it is sufficient to show
{\allowdisplaybreaks
  \begin{align}
    & \sum_{j=0}^{n-1}G_{n,j}(z)\sum_{r=0}^{j}z^{j-r}
      \sum_{l=1}^{r+1}\frac{(-1)^{l-1}z^{l}}{l}k^{r-l}
      \label{eqn;326}\\[8pt]
    & \quad = \frac{1}{(n-1)!}\sum_{\mu = -1}^{n-2}
      \left\{ \sum_{r=\mu+1}^{n-1} 
      \frac{\sideset{_{n-1}}{_r}S}{r-\mu}z^{r-\mu}\right\}
      (-1)^{\mu+1}k^{\mu}.\nonumber
  \end{align}}

Since (\ref{eqn;324}) and
{\allowdisplaybreaks
  \begin{align*}
    &  (z+k)^{j}\log \left(1+\frac{z}{k}\right)
       + \sum_{r=0}^{j}\binom{j}{r}z^{j-r}
       \sum_{l=1}^{r+1}\frac{(-1)^{l}z^{l}}{l} k^{r-l}\\[4pt]
    & \qquad = O(k^{-2}) \quad \mbox{as} \quad k \to \infty,
  \end{align*}}
we have 
{\allowdisplaybreaks
  \begin{align}
    & - \binom{-k}{n-1}\log \left(1+ \frac{z}{k}\right)
      + \sum_{j=0}^{n-1}G_{n,j}(z)\sum_{r=0}^{j}
      \binom{j}{r}z^{j-r}
      \sum_{l=1}^{r+1}\frac{(-1)^{l-1}z^{l}}{l}
      k^{r-l}\label{eqn;327}\\[4pt]
    & \qquad = O(k^{-2}) \quad \mbox{as} \quad k\to \infty,
      \nonumber
  \end{align}}
while we obtain
{\allowdisplaybreaks
  \begin{align}
     &  - \binom{-k}{n-1}\log\left(1+\frac{z}{k}\right)
       + \frac{1}{(n-1)!}\sum_{\mu=-1}^{n-2}\left\{
       \sum_{r=\mu+1}^{n-1}\frac{\sideset{_{n-1}}{_r}S}{r-\mu}
       z^{r-\mu}\right\}(-1)^{\mu+1}k^{\mu}\label{eqn;328}\\[8pt]
     & \quad = - \binom{-k}{n-1}\log \left(1+\frac{z}{k}\right)
       + \frac{1}{(n-1)!}\sum_{r=0}^{n-1}(-1)^{r}
       \sideset{_{n-1}}{_r}S k^{r}\sum_{l=1}^{r+1}
       \frac{(-1)^{l-1}}{l} \left(\frac{z}{k}\right)^{l}
       \nonumber\\[4pt]
     & \quad = O(k^{-2}) \quad \mbox{as} \quad k \to \infty,
       \nonumber
  \end{align}}
by Lemma \ref{lem;34} (1).
We can deduce that (\ref{eqn;327}) and (\ref{eqn;328}) imply (\ref{eqn;326})
by the same arguments as in Lemma 5.8. Hence the proof is completed.
\qed 


\subsection{Examples of the Weierstrass product representation for
$G_{n}(z+1)$}
We give some examples of the Weierstrass product representation for the 
multiple gamma functions.\par
In the case that $n=1$, we have 
  \begin{equation*}
      G_{1}(z+1) = \Gamma(z+1)
        = e^{-\gamma z}\prod_{k=1}^{\infty}
        \left\{\left(1+\frac{z}{k}\right)^{-1}
        e^{-\frac{z}{k}}\right\}.
   \end{equation*}  
This is the Weierstrass product representation for the gamma function.\par
In the case that $n=2$, we have
  \begin{equation*}
    G_{2}(z+1) = G(z+1)
      = e^{-z\zeta'(0)-\frac{z^{2}}{2}\gamma
      - \frac{z^{2}+z}{2}}
      \prod_{k=1}^{\infty} \left\{\left(
        1+\frac{z}{k}\right)^{k}
        \exp\left(-z+\frac{z^{2}}{2k}\right)\right\}.
  \end{equation*}
Since $\zeta'(0)=-\frac{1}{2}\log(2\pi)$, this is the Weierstrass product 
representation for the Barnes $G$-function \cite{bar1} \par 
In the case that $n=3$, $4$ and $5$, we obtain the following results.

\begin{prop}
The Weierstrass product representations in the case that  $n=3$, $4$ and $5$
are as follows: 
{\allowdisplaybreaks
  \begin{align*}
    & G_{3}(z+1)\\[4pt]
    & \quad = \exp\left\{-\frac{z^{3}}{4}+\frac{z^{2}}{8}
      +\frac{7}{24}z + \zeta'(-1)-\frac{z(z-1)}{2}\zeta'(0)
      -\left(\frac{z^{3}}{6}-\frac{z^{2}}{4}\gamma\right)\right\}\\[8pt]
    & \quad \qquad \times \prod_{k=1}^{\infty}\left[
      \left(1+\frac{z}{k}\right)^{-\frac{k(k+1)}{2}}
      \exp\left\{\left(\frac{z^{3}}{6}-\frac{z^{2}}{4}\right)\frac{1}{k}
      -\left(\frac{z^{2}}{4}-\frac{z}{2}\right)
      +\frac{z}{2}k\right\}\right],
\\[12pt]
   & G_{4}(z+1)\\[4pt]
   & \quad = \exp\left\{\frac{61}{144}z^{4}+\frac{13}{18}z^{3}
     + \frac{19}{144}z^{2}-\frac{5}{24}z \right.\\[8pt]
   & \quad \qquad \left.  
     - \frac{z}{2}\zeta'(-2) + \frac{z^{2}-2z}{3}\zeta'(-1)
     - \frac{z^{3}-3z^{2}+2z}{6}\zeta'(0)
     - \frac{z^{4}-4z^{3}+4z^{2}}{24}\gamma\right\}\\[8pt]
   & \quad \qquad\times \prod_{k=1}^{\infty}\left[
     \left(1+\frac{z}{k}\right)^{\frac{k(k+1)(k+2)}{6}}
     \exp\left\{\left(
     \frac{z^{4}}{24}-\frac{z^{3}}{6}+\frac{z^{2}}{6}\right)
     \frac{1}{k}\right.\right.\\[8pt]
   & \quad \qquad \qquad \qquad \qquad \left.\left.
     - \left(\frac{z^{3}}{18}-\frac{z^{2}}{4}-\frac{z}{3}\right)
     + \left(\frac{z^{2}}{12}-\frac{z}{2}\right)k
     - \frac{z}{6}k^{2}\right\}\right],
\\[12pt]
  & G_{5}(z+1)\\[4pt]
  & \quad = \exp\left\{-\frac{5}{288}z^{5}+\frac{7}{64}z^{4}
    -\frac{173}{864}z^{3}-\frac{z^{2}}{36}+\frac{2827}{17280}z\right.\\
  & \quad \qquad + \frac{z}{6}\zeta'(-3)-\frac{z^{2}-3z}{4}\zeta'(-2)
    + \frac{2z^{3}-9z^{2}+11z}{12}\zeta'(-1)\\[8pt]
  & \quad \qquad \left.- \frac{z^{4}-6z^{3}+11z^{2}-6z}{24}\zeta'(0)
    - \frac{6z^{5}-45z^{4}+110z^{3}-90z^{2}}{720}\gamma\right\}\\[8pt]
  & \quad \qquad \times \prod_{k=1}^{\infty}
    \left[\left(1+\frac{z}{k}\right)^{-\frac{k(k+1)(k+2)(k+3)}{24}}
    \exp\left\{\left(\frac{z^{5}}{120}-\frac{z^{4}}{16}+\frac{11}{72}z^{3}
    -\frac{z^{2}}{8}\right)\frac{1}{k}\right.\right.\\[8pt]
  & \quad \qquad \qquad - \left(\frac{z^{4}}{96}-\frac{z^{3}}{12}
    +\frac{11}{48}z^{2}-\frac{z}{4}\right)+\left(\frac{z^{3}}{72}
    -\frac{z^{2}}{8}+\frac{11}{24}\right)k\\[8pt]
  & \quad \qquad \qquad \qquad \left.\left.
    -\left(\frac{z^{2}}{24}-\frac{z}{4}\right)k^{2}
    +\frac{z}{24}k^{3}\right\}\right].
  \end{align*}}  
\end{prop}


\end{document}